\documentclass[12pt,a4paper]{amsart}

\usepackage{color,bm}
\usepackage{amscd,amssymb,amsopn,amsmath,amsthm,graphics,amsfonts,enumerate,verbatim,calc}
\usepackage[dvips]{graphicx}

\makeatletter
  
  \@addtoreset{equation}{section}
\makeatother

\newtheorem{Theorem}{Theorem}[section]
\newtheorem{Proposition}[Theorem]{Proposition}
\newtheorem{Lemma}[Theorem]{Lemma}
\newtheorem{Corollary}[Theorem]{Corollary}
\newtheorem{Remark}[Theorem]{Remark}
\newtheorem{Example}[Theorem]{Example}
\newtheorem{Definition}[Theorem]{Definition}

\newtheorem{Claim}[Theorem]{Claim}

\newcommand{\rmn}[1]{\romannumeral#1}
\def\RMN#1{\uppercase\expandafter{\romannumeral#1}}

\newcommand{\height}{\mathop{\rm ht}\nolimits}

\newcommand{\phideg}{\phi\mathchar`-{\rm deg}}

\newcommand{\oo}{\mathcal{O}}

\newcommand{\bZ}{\mathbb{Z}}

\newcommand{\bQ}{\mathbb{Q}}
\newcommand{\bP}{\mathbb{P}}

\newcommand{\bQp}{\mathbb{Q}_{>0}}

\newcommand{\bN}{\mathbb{N}}
\newcommand{\bNo}{\mathbb{N}_0}
\newcommand{\bR}{\mathbb{R}}
\newcommand{\bC}{\mathbb{C}}
\newcommand{\bRp}{\mathbb{R}_{> 0}}
\newcommand{\bRo}{\mathbb{R}_{\ge 0}}
\newcommand{\bma}{{\bm{a}}}
\newcommand{\bmb}{{\bm{b}}}
\newcommand{\bmd}{{\bm{d}}}
\newcommand{\bme}{{\bm{e}}}
\newcommand{\bmc}{{\bm{c}}}
\newcommand{\bmx}{{\bm{x}}}
\newcommand{\bmy}{{\bm{y}}}
\newcommand{\bmz}{{\bm{z}}}
\newcommand{\bmf}{{\bm{f}}}
\newcommand{\bmg}{{\bm{g}}}

\newcommand{\bmo}{{\bm{0}}}

\begin{document}
\title[Demazure construction]{Demazure construction \\
for $\bZ^n$-graded Krull domains}
\author{Yusuke Arai, Ayaka Echizenya and Kazuhiko Kurano}
\dedicatory{Dedicated to Professor~Kei-ichi~Watanabe on the~occasion of his~74th~birthday.}
\date{}
\maketitle

\begin{abstract}
For a Mori dream space $X$, the Cox ring ${\rm Cox}(X)$ is a Noetherian $\bZ^n$-graded normal domain for some $n > 0$.
Let  $C({\rm Cox}(X))$ be the cone (in $\bR^n$) which is spanned by the vectors $\bma \in \bZ^n$
such that ${\rm Cox}(X)_\bma \neq 0$.
Then  $C({\rm Cox}(X))$ is decomposed into a union of chambers.
Berchtold and Hausen~\cite{BH} proved the existence of such decompositions for 
affine integral domains over an algebraically closed field.
We shall give an elementary algebraic proof to this result
in the case where the homogeneous component of degree $\bmo$ is a field.

Using such decompositions, we develop the Demazure construction for 
$\bZ^{n}$-graded Krull domains.
That is, under an assumption, we show that a $\bZ^n$-graded Krull domain is isomorphic to the multi-section ring $R(X; D_1, \ldots, D_n)$ 
for certain normal projective variety $X$ and $\bQ$-divisors 
$D_1$, \ldots, $D_n$ on $X$.
\end{abstract}

\section{Introduction}

Let $A = \oplus_{m \ge 0}A_{m}$ be a Noetherian normal graded domain over a field $k = A_{0}$.
Then, by the {\em Demazure construction} \cite{D}, \cite{Do}, \cite{P}, \cite{W}~\footnote{Some people call it the {\em Dolgachev-Pinkham-Demazure construction}.},
there exists a normal projective variety $X$ over $k$ and a $\bQ$-divisor $D$ such that
$dD$ is an ample Cartier divisor on $X$ for some $d > 0$, and
$A$ is isomorphic to 
\[
\bigoplus_{r \ge 0} H^{0}(X, \oo_{X}(rD))
\]
as a graded ring.
When $X = {\rm Spec}(k)$, we think that $0$ is the unique $\bQ$-divisor, and 
 it is an ample Cartier divisor.

One of our aims in this paper is to develop the Demazure construction for 
$\bZ^{n}$-graded Krull domains (Theorem~\ref{Demazure} and \ref{Demazure2}).
For a normal projective variety $X$ and $\bQ$-divisors $D_{1}$, \ldots, $D_{n}$, we put
\[
R(X;D_{1}, \ldots, D_{n}) = \bigoplus_{r_{1}, \ldots, r_{n} \in \bZ} H^{0}(X, \oo_{X}(\sum_{i = 1}^{n}r_{i}D_i)) .
\]
It has a $\bZ^{n}$-graded ring structure, and is called the {\em multi-section ring}
with respect to $X$ and $D_{1}$, \ldots, $D_{n}$.
Under an assumption, we show that a $\bZ^n$-graded Krull domain is isomorphic to the multi-section ring $R(X; D_1, \ldots, D_n)$ 
for certain normal projective variety $X$ and $\bQ$-divisors 
$D_1$, \ldots, $D_n$ on $X$.
In the case of $\bZ^n$-graded affine integral domain over an algebraically closed field of characteristic $0$, Altmann and Hausen developed another generalization
of the Demazure construction using proper polyhedral divisors in Theorem~3.4 in \cite{AH}.
In their construction, the given ring itself is not equal to a multi-section ring.
However, using proper polyhedral divisors, they patch together multi-section
rings to obtain the given ring. 

Even if $D_{1}$, \ldots, $D_{n}$ are ample Cartier divisors,
$R(X;D_{1}, \ldots, D_{n})$ is not necessary a Noetherian ring.
We remark that $R(X;D_{1}, \ldots, D_{n})$ is a Noetherin ring
if and only if it is finitely generated over $R(X;D_{1}, \ldots, D_{n})_\bmo$ as a ring.

Finite generation of $R(X;D_{1}, \ldots, D_{n})$ is deeply related to birational geometry.
Let $X$ be a normal $\bQ$-factorial projective variety over a field $k$,
and assume that $D_{1}$, \ldots, $D_{n}$ are free bases of 
the divisor class group ${\rm Cl}(X)$.
Then the ring $R(X;D_{1}, \ldots, D_{n})$ is called the {\em Cox ring} of $X$,
and denoted by ${\rm Cox}(X)$.
We say that $X$ is a MDS (Mori dream space) if ${\rm Cox}(X)$ is Noetherian~\cite{HK}.
Finite generation of Cox rings is a very important problem both in algebraic geometry and in commutative ring theory.
Let $X$ be a blow-up of $\bP^2_\bC$ at finite closed points $p_{1}$, \ldots, $p_{s}$.
Nagata~\cite{Nagata} proved that ${\rm Cox}(X)$ is the invariant subring
of a polynomial ring over $\bC$ with some linear action under some weak assumption.
It is known  that ${\rm Cox}(X)$ is not finitely generated if $s \ge 9$ and 
points $p_{1}$, \ldots, $p_{s}$ are very general.
Therefore it gives a counterexample to Hilbert's 14th problem.

For a MDS $X$, consider the $\bZ^n$-graded ring  ${\rm Cox}(X)$.
We define the cone 
\[
C( {\rm Cox}(X)) = \sum_{\bma \in \bZ^n, \  {\rm Cox}(X)_\bma \neq 0} \bRo   \bma
\subset \bR^n .
\]
The cone $C( {\rm Cox}(X))$ is divided into some chambers and each chamber corresponds to a projective variety 
which is birational to $X$ (cf.\ Hu-Keel~\cite{HK}, Okawa~\cite{Okawa}, Laface-Velasco~\cite{LV}).
In order to define a chamber,  Laface-Velasco studied the ideal generated by elements with degree in a given ray, 
i.e., for $\bma \in \bZ^n$,
\[
J_\bma ({\rm Cox}(X)) = \sqrt{\mbox{the ideal of ${\rm Cox}(X)$ generated by $\cup_{r > 0} {\rm Cox}(X)_{r\bma}$}} .
\]
For $\bma \in \bZ^n$, we put
\[
X_\bma = {\rm Proj}(\oplus_{r \ge 0} {\rm Cox}(X)_{r\bma}) .
\]
It is easy to prove that, if 
\[
J_\bmb ({\rm Cox}(X)) \supset J_\bma ({\rm Cox}(X)),
\]
then we have a morphism
\[
X_\bma \longrightarrow X_\bmb
\]
as in Lemma~\ref{morph}.
We obtain important birational morphisms analyzing ideals of the form $J_\bma ({\rm Cox}(X))$.
One of our aims is to study such ideals for (not necessary Noetherian) $\bZ^n$-graded rings.

In section~2, for a $\bZ^n$-graded ring $A$, 
we define the cone $C(A)$, the ray ideal $J_\bma(A)$ and ray ideal cones which are the generalization of chambers.
We study basic properties of them (cf.\ Proposition~\ref{reduce}).
In Example~\ref{2.2}, we know that $C(A)$ is a union of finitely many chambers 
for a Noetherian $\bZ^2$-graded domain $A$ such that $A_\bmo$ is a field.

In section~3, for a Noetherian $\bZ^n$-graded domain $A$ such that $A_\bmo$ is a field, we shall prove that 
if  $J_\bma (A) = J_\bmb (A)$, then $J_\bmc(A) = J_\bma(A)$ for any $\bmc$ on the line segment between $\bma$ and $\bmb$ (cf.\ Theorem~\ref{eti}
and Corollary~\ref{open}).
This result follows from Theorem~2.11 in \cite{BH}.
If we remove the assumption that $A$ is a Noetherian, Theorem~\ref{eti} is false.
We give a counter example which is a Cox ring of a normal projective rational surface in Example~\ref{SpacceMonomial}.

In section~4, we study Noetherian $\bZ^n$-graded domains with only one chamber.
For a $\bNo^n$-graded Noetherian ring $A$, we give a necessary and sufficient conditions for $A$ to be integral over the subring generated by elements with degree in the coordinate axes
(Theorem~\ref{Arai}). 

In section~5, for a Noetherian $\bZ^n$-graded domain $A$ such that $A_\bmo$ is a field, we decompose $C(A)$
into a union of chambers (cf.\ Theorem~\ref{ChamberDecomp}).

In section~6, we refine arguments in the previous section. 
Considering maximal ray ideal cones in stead of chambers, we give a structure of a fan to the set of maximal ray ideal cones (cf.\ Theorem~\ref{FanStructure}).
It is already known by Theorem~2.11 in \cite{BH}.
For the non-Noetherian symbolic Rees ring in Example~\ref{SpacceMonomial},
the set of maximal ray ideal cones do not form a fan.

Theorems~\ref{eti}, \ref{ChamberDecomp}, \ref{FanStructure} follow from
Theorem~2.11 in \cite{BH}.
However, we give proofs since we need these arguments in the later sections
and our proofs are elementary and algebraic.

 In section~7, we study basic properties of $R(X;D_{1}, \ldots, D_{n})$
for a normal projective variety $X$ and $\bQ$-divisors $D_{1}$, \ldots, $D_{n}$
(cf.\ Theorem~\ref{EKW}).
It is a generalization of results in \cite{EKW}, \cite{K31}, \cite{W}.
We shall prove that $R(X;D_{1}, \ldots, D_{n})$ is a Krull domain and determine the
divisor class group of it.

Using these results, we study the Demazure construction for $\bZ^n$-graded 
Krull domains in Section~8 (cf.\ Theorems~\ref{Demazure} and \ref{Demazure2}).

\section{Notation and basic properties}

Throughout of this paper, $\bC$, $\bR$, $\bQ$, $\bZ$, $\bNo$, and $\bN$ denote 
the set of complex numbers, real numbers, rational numbers, integers,  non-negative integers, and positive integers, respectively.
Furthermore, $\bRo$ and $\bRp$ denote the set of non-negative real numbers and
positive real numbers, respectively.

Let 
\[
A = \bigoplus_{\bma \in \bZ^n}A_{\bma}
\]
be a $\bZ^n$-graded ring, that is, 
each $A_{\bma}$ is an additive subgroup such that $A_{\bma}A_{\bmb}\subset A_{\bma + \bmb}$ for $\bma, \bmb \in \bZ^{n}$.
For $\bma \in \bR^n$, $|\bma|$ denotes the length of the vector $\bma$.
For $0 \neq x \in A_\bma$, we say that the degree of $x$ is $\bma$,
and denote $\deg(x) = \bma$.
We put
\[
C(A) = \sum_{A_{\bma} \neq 0}\bRo   \bma \subset \bZ^n \otimes_{\bZ} \bR = \bR^n .
\]

A subset $\sigma$ in $\bR^{n}$ is called a {\em cone} if the following two conditions are satisfied;
(i) $\bma + \bmb \in \sigma$ if $\bma, \bmb \in \sigma$,
(ii) $\bRo   \bma \subset \sigma$ if $\bma \in \sigma$.

For a cone $\sigma$, $\dim \sigma$ denotes the dimension of $\sigma - \sigma$ as an $\bR$-vector space.
We say that $\sigma$ is a {\em rational polyhedral cone} if
$\sigma$ is spanned by finitely many elements in $\bQ^{n}$.

If $A$ is a finitely generated $\bZ^{n}$-graded ring over $A_{\bmo}$, then $C(A)$ is a rational polyhedral cone.

For $\bma \in \bQ^n$, we define $I_\bma (A)$ to be the ideal of $A$ 
generated by
\[
\bigcup_{\bmb \in \bRp   \bma \cap \bZ^n}A_\bmb .
\]
We put
\[
J_\bma (A) = \sqrt{I_\bma (A)} ,
\]
and call it the {\em ray ideal} of $A$ at $\bma$.
In the case where $A$ is a $\bZ^n$-graded domain, for $\bma \in \bQ^{n}$, 
$I_{\bma}(A) \neq 0$ if and only if $\bma \in C(A)$.

Let $\sigma$ be a cone in $\bR^{n}$ such that $\sigma - \sigma$ is an $\bR$-vector subspace spanned by finitely many elements in $\bQ^{n}$.
We say that $\sigma$ is a {\em ray ideal cone} of $A$
if $J_\bma (A)=J_\bmb (A) \neq 0$ for any $\bma, \bmb \in {\rm rel.int}(\sigma) \cap \bQ^n$, where $ {\rm rel.int}(\sigma)$ denotes the relative interior of $\sigma$.
For a ray ideal cone $\sigma$ of $A$,
$J_\sigma(A)$ denotes $J_\bma(A)$ for $\bma \in {\rm rel.int}(\sigma) \cap \bQ^n$,
and call it the {\em ray ideal} of the ray ideal cone $\sigma$.
We sometimes denote $J_\sigma(A)$ simply by $J_\sigma$
if no confusion is possible.
A ray ideal cone  $\sigma$ is called a {\em chamber} of $A$
if $\dim \sigma = \dim C(A)$.

Let $T$ be an additive  subsemigroup of $\bR^n$ containing $\bmo$.
For example, $T$ is a subgroup of $\bZ^{n}$ or a cone in $\bR^{n}$.
We put
\[
A_T = \bigoplus_{\bma \in T \cap \bZ^n} A_\bma .
\]
Remark that it is a subring of $A$.
We regard $A_{T}$ as a $\bZ^{n}$-graded subring of $A$
unless otherwise specified.

\begin{Example}\label{2dim-chamber}
\begin{rm}
Let $A=k[x,y,z]$ be a $\bZ^{2}$-graded polynomial ring 
over a field $k$ with
$\deg(x) = (1,0)$,
$\deg(y) = (1,1)$,
$\deg(z) = (0,1)$,
and $A_{\bmo} = k$.
Then the set of non-zero ray ideals of $A$ consists of
\[
A, \ (x), \ (xz,y), \ (z), \ (xy, xz), \ (xz, yz) .
\]
The following are the maximal ray ideal cones of the above ray ideals respectively.
\[
\{ 0 \}, \ \bRo   (1,0), \ \bRo   (1,1), \ \bRo   (0,1), \ \bRo   (1,0)+\bRo   (1,1), \ \bRo   (1,1)+\bRo   (0,1)
\]
\end{rm}
\end{Example}

\begin{Example}\label{2.2}
\begin{rm}
Let $A = k[x_{1}, x_{2}, \ldots, x_{t}]$ be a $\bZ^{2}$-graded domain
over a field $A_{\bmo} = k$, where $x_{i}$ is a non-zero homogeneous element with
$\deg(x_{i}) = (\alpha_{i},\beta_{i})$ for $i = 1, \ldots, t$.
Here, we do not have to assume that $x_{1}$, $x_{2}$, \ldots, $x_{t}$ are algebraically independent over $k$.
Then $C(A)$ is the cone spanned by $\{ (\alpha_i, \beta_i) \mid i = 1, \ldots, t \}$.
For $\bma \in \bQ^n$, $J_\bma(A) \neq 0$ if and only if $\bma \in C(A)$
since $A$ is a domain.

Consider the following three cases:
\begin{itemize}
\item[(\RMN{1})]
Assume that $C(A)$ is {\em strongly convex}, that is, $-\bma \not\in C(A)$
for any $\bmo \neq \bma \in C(A)$.
Changing coordinates, we may assume $\alpha_{i} > 0$ for all $i$ and
\[
- \infty < \frac{\beta_{1}}{\alpha_{1}} \le \frac{\beta_{2}}{\alpha_{2}} \le \cdots \le
\frac{\beta_{t}}{\alpha_{t}} < \infty .
\]
Assume 
\[
\frac{\beta_{i}}{\alpha_{i}} < \frac{\beta_{i+1}}{\alpha_{i+1}} =  \cdots = \frac{\beta_{i+s}}{\alpha_{i+s}} < \frac{\beta_{i+s+1}}{\alpha_{i+s+1}} .
\]
Then $\bRo   (\alpha_{i},\beta_{i}) + \bRo   (\alpha_{i+1},\beta_{i+1})$ is a chamber
with ray ideal 
\begin{equation}\label{ideal1}
\sqrt{(x_{1}, x_{2}, \ldots, x_{i}) \cap (x_{i+1}, x_{i+2}, \ldots, x_{t})} .
\end{equation}
Furthermore, $\bRo   (\alpha_{i+1},\beta_{i+1})$ is a ray ideal cone with ray ideal 
\begin{equation}\label{ideal2}
\sqrt{(x_{1}, x_{2}, \ldots, x_{i+s}) \cap (x_{i+1}, x_{i+2}, \ldots, x_{t})} .
\end{equation}
Here, remark that the ideal (\ref{ideal1}) is contained in the ideal (\ref{ideal2}).
If
\[
\frac{\beta_{1}}{\alpha_{1}} = \frac{\beta_{2}}{\alpha_{2}}
= \cdots = \frac{\beta_{i}}{\alpha_{i}} < \frac{\beta_{i+1}}{\alpha_{i+1}} ,
\]
then $\bRo   (\alpha_{1},\beta_{1})$ is a ray ideal cone with ray ideal 
\begin{equation}\label{ideal3}
\sqrt{(x_{1}, x_{2}, \ldots, x_{i})} .
\end{equation}
Thus, in this case, $C(A)$ is a union of some chambers $\sigma_{1}$, $\sigma_{2}$, \ldots, $\sigma_{s}$.
Furthermore, for $i = 1, 2, \ldots, s$, the ray ideal of $\sigma_{i}$ is contained in the ray ideals of the faces
of $\sigma_{i}$.
\item[(\RMN{2})]
Assume that $C(A) = \bRo \times \bR$.
Put $\bma = (\alpha, \beta) \in \bQ^2$.
If $\alpha > 0$, then
\[
J_{\bma}(A) = J_{(1,0)}(A) \subsetneq A .
\]
For any $\beta \in \bQ$, 
\[
J_{(0,\beta)}(A) = A .
\]
If $\alpha < 0$, then $J_{\bma}(A) = 0$.
In this case, $A$ has two maximal ray ideal cones $\bRo \times \bR$ and $\{ 0 \} \times \bR$.
\item[(\RMN{3})]
Assume that $C(A) = \bR^2$.
For any $\bma \in \bQ^2$, $J_{\bma}(A) = A$ in this case.
\end{itemize}
\end{rm}
\end{Example}

\begin{Remark}\label{Noetherian}
\begin{rm}
Let 
\[
A = \bigoplus_{\bma \in \bZ^n}A_{\bma}
\]
be a $\bZ^n$-graded ring.
\begin{enumerate}
\item
It is well known that $A$ is Noetherian
if and only if 
\begin{itemize}
\item
$A_\bmo$ is Noetherian, and
\item
$A$ is finitely generated over $A_\bmo$ as a ring.
\end{itemize}
\item
Let $\bmb \in \bQ^n$.
If $A$ is Noetherian, then so is  
\[
\bigoplus_{(\bma,\bmb) \ge 0} A_\bma ,
\]
where $(\bma,\bmb)$ is the inner product of $\bma$ and $\bmb$.
Using it, we can prove that 
if $A$ is Noetherian, then so is $A_\sigma$ for any rational polyhedral cone
$\sigma$ in $\bR^n$,
where $A_\sigma = \oplus_{\bma \in \sigma \cap \bZ^n}A_\bma$.
\item
Suppose that $T$ is a subgroup or a finitely generated sub-monoid of $\bZ^n$.
If $A$ is Noetherian, then so is  $A_T$.
\end{enumerate}
\end{rm}
\end{Remark}

\begin{Proposition}\label{reduce}
Let $A$ be a $\bZ^n$-graded ring.
\begin{enumerate}
\item
Suppose that $S$ is a  $\bZ^n$-graded subring of $A$.
Let $\bma, \bmb \in \bQ^n$.
Assume that $J_\bma(S) \subset J_\bmb(S)$, and
$\sqrt{I_\bma(S) A} = J_\bma(A)$.
Then $J_\bma(A) \subset J_\bmb(A)$ holds.
\item
Let $T$ be a subgroup of $\bZ^n$.
Let $\bma, \bmb \in T$.
Then
\[
J_\bma(A_T) \subset J_\bmb(A_T)
\Longleftrightarrow
J_\bma(A) \subset J_\bmb(A) .
\]
\item
Let $\sigma$ be a cone in $\bR^{n}$.
Then, for $\bma \in {\rm rel.int}(\sigma) \cap \bQ^n$ and $\bmb \in \sigma \cap \bQ^n$,
\[
J_\bma(A_\sigma) \subset J_\bmb(A_\sigma)
\Longleftrightarrow
J_\bma(A) \subset J_\bmb(A) .
\]
In particular, if $\sigma$ is a ray ideal cone of $A$,  
then $\sigma$ itself is a ray ideal cone of $A_\sigma$.
\item
Let $S$ be a  $\bZ^n$-graded subring of $A$.
Assume that $A$ is integral over $S$.
Then, for any $\bma, \bmb \in \bQ^n$,
\[
J_\bma(S) \subset J_\bmb(S)
\Longleftrightarrow
J_\bma(A) \subset J_\bmb(A) .
\]
\end{enumerate}
\end{Proposition}

\proof
%We shall prove (1).
%It is enough to show $I_\bma(S) \subset J_\bmb(A)$.
%For $x \in I_\bma(S)$, there exists $m \in \bN$ such that $x^m \in I_\bmb(S)$
%since $I_\bma(S) \subset J_\bmb(S) = \sqrt{I_{\bmb}(S)}$.
%Then we have $x^m \in I_\bmb(S) \subset I_\bmb(A) \subset J_\bmb(A)$.
The assertion~(1) follows from
\[
J_\bma(A) = \sqrt{I_\bma(S)A} = \sqrt{J_\bma(S)A} \subset \sqrt{J_\bmb(S)A}
\subset J_\bmb(A) . 
\]

We shall prove (2).
Assume $\bma, \bmb \in T$.
For a homogeneous element $w$ of $A$ with degree on the ray $\bRp   \bma$,
some power of $w$ is contained in $I_\bma(A_T)$.
Therefore we have $\sqrt{I_\bma(A_T)A} = J_\bma(A)$.
The implication $(\Rightarrow)$ follows from (1) as above.
Next, we shall prove  $(\Leftarrow)$.
For $x \in A_\bmc$ where $\bmc \in \bRp   \bma \cap T$,
we want to show $x^m \in I_\bmb(A_T)$ for some $m$.
Since $x \in A_{\bmc} \subset J_\bma(A) \subset J_\bmb(A) = \sqrt{I_\bmb(A_T)A}$,
there exists homogeneous elements $y_1$, \ldots, $y_t$, $z_1$, \ldots, $z_t$ of $A$
such that $x^m = \sum_i y_iz_i$ for some $m$, where $\deg(y_1)$, \ldots, $\deg(y_t)$ are in $\bRp   \bmb \cap T$.
We may assume that $m \bmc = \deg(y_i)+\deg(z_i)$ for each $i$.
Then, since $\deg(z_i) = m \bmc - \deg(y_i) \in T$,
$z_{i}$ is in $A_{T}$ for each $i$.
Thus $x^m$ is in $I_\bmb(A_T)$.

We shall prove (3).
The implication $(\Rightarrow)$ follows from (1) as above.
Next, we shall prove  $(\Leftarrow)$.
It is enough to show $A_{\bmc} \subset J_{\bmb}(A_{\sigma})$ for $\bmc \in \bRp   \bma \cap \bZ^{n}$.
Let $x$ be an element of $A_{\bmc}$.
By the assumption, $x \in J_{\bma}(A) \subset J_{\bmb}(A) = \sqrt{I_{\bmb}(A)}$.
Then there exists homogeneous elements $y_1$, \ldots, $y_t$, $z_1$, \ldots, $z_t$ of $A$
such that $x^m = \sum_i y_iz_i$ for some $m$, where $\deg(y_1)$, \ldots, $\deg(y_t)$ are in $\bRp   \bmb$.
We may assume that $m \bmc = \deg(y_i)+\deg(z_i)$ for each $i$.
Let $s$ be a positive integer.
We have
\[
x^{m+s} = \sum_i y_i(x^{s}z_i) .
\]
Since $\deg(z_i) = m \bmc - \deg(y_i)$,
$\deg(z_{i})$ is contained in $\sigma - \sigma$.
Since $\deg(x)$ is in ${\rm rel. int}(\sigma)$,
$\deg(x^{s}z_{i})$ is in $\sigma$ for $s \gg 0$.
Since $x^{m+s} \in I_{\bmb}(A_{\sigma})$ for $s \gg 0$,
we have $x \in J_{\bmb}(A_{\sigma})$.

We shall prove (4).
It is enough to show that
\begin{itemize}
\item[(i)]
$\sqrt{I_{\bma}(S) A} = J_{\bma}(A)$, and
\item[(ii)]
$J_{\bma}(S) = J_{\bma}(A) \cap S$
\end{itemize}
for any $\bma \in \bQ^{n}$.

First, we shall prove (i).
It is easy to see that the left-hand-side is contained in the right one.
In order to show the opposite inclusion,
we shall prove $A_{\bmc} \subset \sqrt{I_{\bma}(S)A}$ for $\bmc \in \bRp   \bma \cap \bZ^{n}$.
Take $x \in A_{\bmc}$.
Since $A$ is integral over $S$,
we have an integral equation
\[
x^{m} + a_{1}x^{m-1} + \cdots + a_{m} = 0
\]
for some $m$ and some $a_{1}, a_{2}, \ldots, a_{m} \in S$.
Here we may assume that $a_{i}$ is contained in $S_{i\bmc}$ for each $i$.
Therefore $x^{m}$ is in $I_{\bma}(S)A$.

Next, we shall prove (ii).
It is easy to see that the left-hand-side is contained in the right one.
We shall prove the opposite inclusion.
It suffices to show that, if a prime ideal $P$ of $S$
contains the left-hand-side, then $P$ contains the right one.
Since $S \hookrightarrow A$ is an integral extension,
there exists a prime ideal $Q$ of $A$ such that $Q \cap S = P$.
Take $x \in A_{\bmc}$ 
for some $\bmc \in \bRp   \bma \cap \bZ^{n}$.
Take an integral equation
\[
x^{m} + a_{1}x^{m-1} + \cdots + a_{m} = 0
\]
for some $m$ and some $a_{1}, a_{2}, \ldots, a_{m} \in S$.
We may assume that $a_{i}$ is contained in $S_{i\bmc}$ for each $i$.
Then $a_{i} \in J_{\bma}(S) \subset P \subset Q$.
By the above integral equation, we know $x \in Q$.
Thus $Q$ contains $J_{\bma}(A)$.
Hence $P$ contains $J_{\bma}(A) \cap S$.
\qed

\begin{Example}\label{xy-zw}
\begin{rm}
Let $A = k[x,y,z,w]/(xw-yz)$ be a $\bZ^{2}$-graded ring with $\deg(x) = (3,0)$, $\deg(y) = (2,1)$, $\deg(z) = (1,2)$, $\deg(w) = (0,3)$ and $A_{\bmo} = k$.
As in Example~\ref{2dim-chamber}, we have mutually distinct 
non-zero ray ideals of $A$ as 
\begin{align*}
& A, \\
& \sqrt{(x)}, \ \sqrt{(y, xz, xw)}, \ \sqrt{(z, xw, yw)}, \ \sqrt{(w)}, \\
& \sqrt{(x) \cap (y,z,w)}, \ \sqrt{(x,y) \cap (z,w)}, \ \sqrt{(x,y,z) \cap (w)} .
\end{align*}

Now, suppose that there exists homogeneous elements $a$, $b$, $c$ in $A$ such that the inclusion
\[
S = k[a,b,c] \rightarrow A
\]
is finite.
Then, by Example~\ref{2.2} (\RMN{1}), $S$ has at most 6 non-zero ray ideals.
Then, by Proposition~\ref{reduce} (4), $A$ also has at most 6 non-zero ray ideals.
It is a contradiction.
Therefore $A$ never have a $\bZ^2$-graded Noether normalization.
\end{rm}
\end{Example}

\section{Ray ideals for $\bZ^{n}$-graded Noetherian rings}

The aim of this section is to prove the following theorem.
This result follows from Theorem~2.11 in \cite{BH}.
If we remove the assumption that $A$ is Noetherian,
it is not true as in Example~\ref{SpacceMonomial}.

\begin{Theorem}\label{eti}
Let $A$ be a Noetherian $\bZ^{n}$-graded domain such that $A_{\bmo}$ is a field~\footnote{
Instead of assuming that $A$ is Noetherian,
it is sufficient to assume that $A_{\bRo}\bma$ is Noetherian.
In fact, if $A_{\bRo}\bma$ is Noetherian,
we can find a Noetherian subring of $A$ with the same situation as $A$.}.
Suppose $\bma, \bmb, \bmc \in \bQ^{n} \subset \bR^{n}$.
Assume that $\bmc$ is on the line segment between $\bma$ and $\bmb$.
Assume $0 \neq J_{\bma}(A) \subset J_{\bmb}(A)$.
Then $J_{\bmc}(A) = J_{\bma}(A)$ if $\bmc \neq \bmb$.
\end{Theorem}

\proof
If $\bmc = \bma$, then it is obvious.
Assume $\bmc \neq \bma, \bmb$ in the rest of this section.

First, we shall prove $I_{\bma}(A) I_{\bmb}(A) \subset J_{\bmc}(A)$.
Let $x \in A_{\bmd}$ and $y \in A_{\bme}$ for $\bmd = u \bma \in \bZ^{n}$ and $\bme = v \bmb \in \bZ^{n}$,
where $u, v \in \bQ_{>0}$.
By the assumption, we have $\bmc = s\bma + (1-s)\bmb$ for a rational number $s$ such that $0 < s < 1$.
Choose an integer $m \ge 0$ such that $ms = uu'$ and $m(1-s) = vv'$ for some $u', v' \in \bN$.
Then $x^{u'}y^{v'} \in A_{m\bmc} \subset I_{\bmc}(A)$.
Thus $xy \in J_{\bmc}(A)$. 

Since $I_{\bma}(A) I_{\bmb}(A) \subset J_{\bmc}(A)$, we have
\[
J_{\bma}(A) \subset J_{\bmc}(A)  .
\]

\vspace{2mm}

In the rest of this section, we shall show the opposite inclusion.

If $\bma$ and $\bmb$ are not linearly independent over $\bR$,
it is easy to prove the assertion.

Assume that $\bma$ and $\bmb$ are linearly independent over $\bR$.
Put $T = (\bR \bma + \bR \bmb) \cap \bZ^{n}$.
Replacing $A$ by $A_{T}$,
we may assume that $A$ is $\bZ^{2}$-graded
by Proposition~\ref{reduce} (2).
Consider the normalization $\tilde{A}$ of $A$.
It is well known that $\tilde{A}$ also has a structure of a $\bZ^{2}$-graded ring.
Replacing $A$ by $\tilde{A}$, we may assume that $A$ is a $\bZ^{2}$-graded
normal domain by Proposition~\ref{reduce} (4).

In this case, $C(A)$ is a $2$-dimensional cone
such that $\bma, \bmb \in C(A) \subset \bR^{2}$.
Let $\sigma$ be a strongly convex cone in $\bR^2$ such that
$\bma, \bmb \in {\rm int}(\sigma)$.
By Proposition~\ref{reduce} (3), we may replace $A$ by $A_\sigma$.
In the rest, we assume that $C(A)$ is strongly convex.

As in Example~\ref{2.2} (I), one of the following two cases occurs:
\begin{itemize}
\item[(i)]
There exists a chamber $\tau$ such that $\bma \in {\rm int}(\tau)$.
\item[(ii)]
There exists a chamber $\tau$ which is contained in 
$\bRo \bma + \bRo \bmb$
such that $\bma$
is in the boundary of $\tau$.
\end{itemize}
Assume that (ii) occurs.
Then $J_{\bma}(A)$ contains $J_{\tau}(A)$
by Example~\ref{2.2}.
Then there exists a rational point $\bma'$ (near $\bma$) on the line segment between $\bma$ and $\bmb$ such that $\bma'$ is in the interior of $\tau$.
Then we have $J_{\bma'}(A) \subset J_{\bma}(A)$.
It is enough to show $J_{\bma'}(A) \supset J_{\bmc}(A)$.
Therefore, replacing $\bma$ by $\bma'$, we may assume that
(i) as above occurs.

In the rest, let $A$ be a $\bZ^{2}$-graded normal domain and we assume that $\bma = (\alpha,0)$ and $\bmb = (\gamma, \delta)$, where 
$\alpha$, $\gamma$, $\delta$ are positive rational numbers.
Put $\bmc = s \bma + (1-s)\bmb = (s\alpha + (1-s)\gamma, (1-s) \delta)$, where $s$ is a rational number such that $0 < s < 1$.
Here, remark that
\begin{equation}\label{bmc}
0 < \frac{ (1-s) \delta}{s\alpha + (1-s)\gamma} < \frac{\delta}{\gamma} .
\end{equation}
Since $\bma$ is in the interior of the chamber $\tau$,
there exists $\alpha_1, \beta_1 \in \bQp$ such that 
\begin{equation}\label{rho}
\left( \bRo   \bma + \bRo  (\alpha_1, \beta_1) \right) \setminus \{ \bmo \}
\subset {\rm int}(\tau) .
\end{equation}
Put
\[
\rho = \bRo   \bma + \bRo  (\alpha_1, \beta_1) .
\]
Since $A_\rho$ is Noetherian, there exists homogeneous elements 
$y_1$, \ldots, $y_t$ such that 
\[
A_{\rho} = A_{\bmo}[y_{1}, y_{2}, \ldots, y_{t}] .
\]
Put $\deg(y_{i}) = (\gamma_{i}, \delta_{i}) \in \rho \setminus \{ \bmo \}$ for $i = 1, 2, \ldots, t$.
By definition, we have $\gamma_{i} > 0$ and $\delta_{i} \ge 0$ for $i = 1, 2, \ldots, t$.
Remark that 
\[
0 \le \frac{\delta_{i}}{\gamma_{i}} \le \frac{\beta_1}{\alpha_1}
\]
for each $i$.

Let $x$ be a homogeneous element of $A$ such that
$\deg(x) \in \bRp   \bmc \cap \bZ^2$.
It is enough to show $x \in J_\bma(A)$.
Since $\tau$ is a chamber,
\[
y_{1}, y_{2}, \ldots, y_{t} \in J_\bma(A) \subset J_\bmb(A) 
\]
by (\ref{rho}).
Therefore there exist homogeneous elements $b_{1}$, $b_{2}$, \ldots, $b_{h}$ such that
\begin{equation}\label{123}
\begin{array}{l}
{\rm (\rmn{1})} \
y_{1}, y_{2}, \ldots, y_{t} \in \sqrt{(b_{1}, b_{2}, \ldots, b_{h})A} ,  \\
{\rm (\rmn{2})} \ 
\mbox{
$\deg(b_{1})= \deg(b_{2})= \cdots = \deg(b_{h}) = e\bmb = (e\gamma, e\delta)$ for some $e \in \bN$,} \\
{\rm (\rmn{3})} \ 
\mbox{
$e\gamma$ is strictly bigger than the first component of the vector $\deg(x)$.}
\end{array}
\end{equation}
Hence $y_{1}^{w}$, $y_{2}^{w}$, \ldots, $y_{t}^{w}$ are in $(b_{1}, b_{2}, \ldots, b_{h})A$ for some $w \in \bN$.
We put
\[
\varphi = \min\{ \delta_{i} - \left( \beta_{1}/\alpha_{1} \right) \gamma_{i} \mid
i = 1, 2, \ldots, t \} 
\]
and
\[
\Omega = \left\{ (\xi,\eta) \in \bR^{2} \ \left| \
0 < \xi, \ 0 \le \eta \le \left( \beta_{1}/\alpha_{1} \right) \xi + w t \varphi \right.
\right\} .
\]
By definition, we know $\varphi < 0$.
Hence we have $\Omega \subset \rho$.

Here, we shall prove the following claim:

\begin{Claim}\label{omega}
If $\bmd \in \Omega \cap \bZ^{2}$, then $A_{\bmd} \subset (b_{1}, b_{2}, \ldots, b_{h})A$.
\end{Claim}

Put $\bmd = (\xi_{1}, \eta_{1}) \in \Omega \cap \bZ^{2}$.
Remark that $A_{\bmd} = (A_{\rho})_{\bmd}$.
Put
\[
M = y_{1}^{u_{1}} y_{2}^{u_{2}} \cdots y_{t}^{u_{t}} \in A_{\bmd} .
\]
Here, $(\xi_{1}, \eta_{1}) = (u_{1}\gamma_{1} + u_{2}\gamma_{2} + \cdots + u_{t}\gamma_{t}, u_{1}\delta_{1} + u_{2}\delta_{2} + \cdots + u_{t}\delta_{t})$.
If $M$ is not in $(b_{1}, b_{2}, \ldots, b_{h})A$, then
$u_{i} < w$ for $i = 1,  2, \ldots, t$.
Remark that 
\[
0 \ge \delta_{i} - \frac{\beta_{1}}{\alpha_{1}} \gamma_{i} \ge \varphi
\]
for each $i$.
Then
\begin{eqnarray*}
\eta_{1} - \frac{\beta_{1}}{\alpha_{1}} \xi_{1} & = & 
(u_{1}\delta_{1} + u_{2}\delta_{2} + \cdots + u_{t}\delta_{t})
- \frac{\beta_{1}}{\alpha_{1}} (u_{1}\gamma_{1} + u_{2}\gamma_{2} + \cdots + u_{t}\gamma_{t}) \\
& = & u_{1}\left( \delta_{1} - \frac{\beta_{1}}{\alpha_{1}} \gamma_{1} \right)
+ u_{2}\left( \delta_{2} - \frac{\beta_{1}}{\alpha_{1}} \gamma_{2} \right)
+ \cdots + u_{t}\left( \delta_{t} - \frac{\beta_{1}}{\alpha_{1}} \gamma_{t} \right) \\
& \ge & (u_1 + u_2 + \cdots + u_t) \varphi > w t \varphi .
\end{eqnarray*}
It contradicts $\bmd \in \Omega$.
We have completed the proof of Claim~\ref{omega}.

\vspace{2mm}

Consider the subring 
\[
A_{\bRo   \bma} = \bigoplus_{\bmg \in \bRo   \bma \cap \bZ^n} A_\bmg
\]
 of $A$.
Remark that $A_{\bRo   \bma}$ is Noetherian by Remark~\ref{Noetherian}.
We think $A_{\bRo   \bma}$ as a $\bZ^{2}$-graded subring of $A$.
Let $z_{1}$, $z_{2}$, \ldots, $z_{p}$ be homogeneous elements
satisfying $\deg(z_i) \in \bRp   \bma \cap \bZ^2$ for $i = 1, 2, \ldots, p$ and
\[
A_{\bRo   \bma} = A_{\bmo}[z_{1}, z_{2}, \ldots, z_{p}] .
\]

Then there exist positive integers $k_{1}$, $k_{2}$, \ldots, $k_{p}$ such that
\[
\begin{array}{l}
{\rm (i)} \ 
\mbox{$\deg(z_{1}^{k_{1}}) = \deg(z_{2}^{k_{2}}) = \cdots = \deg(z_{p}^{k_{p}})$, and} \\
{\rm (ii)} \ 
\mbox{$\deg(xz_{j}^{k_{j}}) \in \Omega$ for $j = 1, 2, \ldots, p$.}
\end{array}
\]

By Claim~\ref{omega}, we have
\[
xz_{j}^{k_{j}} = \sum_{u = 1}^{h} b_{u}c_{ju}
\]
where $c_{ju} = 0$ or $c_{ju}$ is a homogeneous element of $A$ with
$\deg(xz_{j}^{k_{j}}) = \deg(b_{u}c_{ju})$ for each $u$.
Assume $c_{ju} \neq 0$.
Then the second component of $\deg(c_{ju})$ is negative 
by (\ref{bmc}) and the definition of $b_{1}$, $b_{2}$, \ldots, $b_{h}$ (cf.\ (\ref{123}) (iii)).
Let $\zeta$ be the length of the vector $e \bmb = \deg(b_{u})$, which is independent of $u$.
Let $L$ be the line segment between $\deg(xz_{j}^{k_{j}})$ and $\deg(c_{ju})$.
Let $\bmf$ be the intersection of $L$ and the coordinate axis $\bR\bma$.
Let $\nu\zeta$ be the distance between $\bmf$ and $\deg(xz_{j}^{k_{j}})$.
Then $\nu$ is a rational number such that $0 < \nu < 1$.
Put 
\[
\phi = |\bmf - \deg(z_{j}^{k_{j}})| ,
\]
where $| \ |$ denotes the length of a vector.
Let $\mu$ be a positive integer such that
\begin{equation}\label{12'}
\begin{array}{l}
{\rm (i)} \ \mbox{
$\mu\nu$ is a positive integer, and} \\
{\rm (ii)} \ \mbox{
$|\deg(z_{1}^{k_{1}} z_{2}^{k_{2}} \cdots z_{p}^{k_{p}})| < \mu \phi$.}
\end{array}
\end{equation}
Since $\mu > \mu\nu >0$, 
we can describe as 
\[
(xz_{j}^{k_{j}})^{\mu}
= \left( \sum_{u = 1}^{h} b_{u}c_{ju} \right)^{\mu}
= \sum_{w_{1}+w_{2} +\cdots + w_{h} = \mu\nu} b_{1}^{w_{1}}b_{2}^{w_{2}}\cdots b_{h}^{w_{h}} C_{j,w_{1},w_{2},\ldots, w_{h}} ,
\]
where $C_{j,w_{1},w_{2},\ldots, w_{h}}$ is a homogeneous element of $A$ with 
\begin{equation}\label{3.4}
\deg(C_{j,w_{1},w_{2},\ldots, w_{h}}) = \mu \deg(xz_{j}^{k_{j}}) - \mu\nu e \bmb =
\mu\bmf .
\end{equation}
Therefore
\begin{equation}\label{DEGofC}
|\deg(C_{j,w_{1},w_{2},\ldots, w_{h}})| = \mu |\deg(z_{j}^{k_{j}})| + \mu \phi .
\end{equation}

On the other hand, by (\ref{3.4}), $C_{j,w_{1},w_{2},\ldots, w_{h}}$ is an $A_{\bmo}$-linear combination of monomials $z_{1}^{v_{1}}z_{2}^{v_{2}}\cdots z_{p}^{v_{p}}$
with $\deg(z_{1}^{v_{1}}z_{2}^{v_{2}}\cdots z_{p}^{v_{p}}) = \mu\bmf$.
We put $v_{j} = r_{j}k_{j} + q_{j}$ for $r_{j}, q_{j} \in \bZ$ such that $0 \le q_{j} < k_{j}$ for $j = 1, 2, \ldots, p$.
Then we have
\[
|\deg(z_{1}^{q_{1}} z_{2}^{q_{2}} \cdots z_{p}^{q_{p}})| < |\deg(z_{1}^{k_{1}} z_{2}^{k_{2}} \cdots z_{p}^{k_{p}})|
< \mu\phi
\]
by (\ref{12'}).
Therefore
\begin{eqnarray*}
|\deg(C_{j,w_{1},w_{2},\ldots, w_{h}})| & = &  |\deg(z_{1}^{v_{1}}z_{2}^{v_{2}}\cdots z_{p}^{v_{p}})| \\
& = & |\deg(z_{1}^{r_{1}k_{1}}z_{2}^{r_{2}k_{2}}\cdots z_{p}^{r_{p}k_{p}})| + |\deg(z_{1}^{q_{1}}z_{2}^{q_{2}}\cdots z_{p}^{q_{p}})| \\
& < & |\deg( z_{1}^{r_{1}k_{1}}z_{2}^{r_{2}k_{2}}\cdots z_{p}^{r_{p}k_{p}})| + \mu\phi .
\end{eqnarray*}
Comparing the above inequality with the equation (\ref{DEGofC}), we have 
\[
\mu |\deg(z_{j}^{k_{j}})| < 
|\deg( z_{1}^{r_{1}k_{1}}z_{2}^{r_{2}k_{2}}\cdots z_{p}^{r_{p}k_{p}})| =
(r_{1}+r_{2}+ \cdots + r_{p})|\deg(z_{j}^{k_{j}})|.
\]
Therefore we obtain
\[
\mu < r_{1}+r_{2}+ \cdots + r_{p} .
\]
Then
\[
C_{j,w_{1},w_{2},\ldots, w_{h}} \in (z_{1}^{k_{1}}, z_{2}^{k_{2}}, \ldots, z_{p}^{k_{p}})^{\mu+1}A .
\]
Let $V$ be a discrete valuation ring containing $A$.
Let $v$ be the valuation of $V$.
Put 
\[
v_0 = v(z_{j}^{k_{j}}) = \min\{ v(z_{1}^{k_{1}}), v(z_{2}^{k_{2}}), \ldots, v(z_{p}^{k_{p}})\} .
\]
Since $(xz_{j}^{k_{j}})^{\mu} \in (z_{1}^{k_{1}}, z_{2}^{k_{2}}, \ldots, z_{p}^{k_{p}})^{\mu+1}A$,
we know that 
\[
v(x^{\mu}) \ge v_0 .
\]
Therefore $x^{\mu}$ is contained in the integral closure of $(z_{1}^{k_{1}}, z_{2}^{k_{2}}, \ldots, z_{p}^{k_{p}})A$ (cf.\ (6.8.2) in \cite{HS}).
In particular, 
\[
x \in \sqrt{(z_{1}^{k_{1}}, z_{2}^{k_{2}}, \ldots, z_{p}^{k_{p}})A} = J_{\bma}(A) .
\]
We have completed the proof of Theorem~\ref{eti}.
\qed

\begin{Corollary}\label{open}
Let $A$ be a Noetherian $\bZ^n$-graded domain such that
$A_{\bmo}$ is a field.
Let $\bma$ and  $\bmb$ be linearly independent vectors  in $\bQ^n$ such that  
$J_{\bma}(A) = J_{\bmb}(A) \neq 0$.
Consider the line $L$ passing through  $\bma$ and  $\bmb$.

Then there exists a sufficiently small $\epsilon$ such that
$J_{\bmc}(A)$ coincides with $J_{\bmb}(A)$ for any $\bmc \in L \cap U_{\epsilon}(\bmb) \cap \bQ^n$,
where $U_{\epsilon}(\bmb) = \{ \bmd \in \bR^n \mid |\bmd - \bmb| < \epsilon \}$.
\end{Corollary}

\proof
Let $T = (\bR \bma + \bR \bmb) \cap \bZ^2$.
Replacing $A$ by $A_T$, we may assume that 
$A$ is a Noetherian $\bZ^2$-graded domain.
Take $\bmc \in L \cap U_{\epsilon}(\bmb) \cap \bQ^n$
for a sufficiently small $\epsilon$.
If $\bmc$ is on the line segment between  $\bma$ and  $\bmb$,
then $J_{\bmc}(A) = J_{\bmb}(A)$ by Theorem~\ref{eti}.
Assume that  $\bmc$ is not on the line segment between  $\bma$ and  $\bmb$.
As in Example~\ref{2.2}, one of the following three cases occurs:
\begin{itemize}
\item[(i)]
${\bmb}$ is in the interior of a chamber $\tau$.
\item[(ii)]
$\bmb$ is in the boundary of two chambers $\sigma_1$ and $\sigma_2$.
\item[(iii)]
$\bmb$ is in the boundary of $C(A)$.
\end{itemize}
If (i) occurs, then we may assume that $\bmc$ is in the interior of the chamber $\tau$.
Hence $J_\bmb(A) = J_\bmc(A)$.
If (ii) occurs, then 
\[
0 \neq J_\bmc(A) \subset J_\bmb(A) = J_\bma(A) .
\]
By Theorem~\ref{eti}, we have $J_\bmc(A) = J_\bmb(A)$.
Assume that (iii) occurs.
Let $x$ be a non-zero homogeneous element such that $\deg(x) = s \bmb$ for some $s \in \bN$.
Then any power of $x$ is not in $I_{\bma}(A)$.
It contradicts $J_{\bma}(A) = J_{\bmb}(A)$.
\qed

\vspace{2mm}

Theorem~\ref{eti} is false if we remove the assumption that $A$ is Noetherian.

\begin{Example}\label{SpacceMonomial}
\begin{rm}
Let $a$, $b$, $c$ be pairwise coprime positive integers.
Put $S = k[x,y,z]$, where $k$ is a field.
Let $P$ be the kernel of the $k$-algebra homomorphism
\[
\phi : S \longrightarrow k[T]
\]
defined by $\phi(x) = T^{a}$, $\phi(y) = T^{b}$ and $\phi(z) = T^{c}$.
For $n \in \bN$, put 
\[
P^{(n)} = P^{n}S_{P} \cap S
\]
and call it the {\em $n$th symbolic power} of $P$.
Put
\[
A = S[t^{-1}, Pt, P^{(2)}t^{2}, P^{(3)}t^{3}, \ldots ] \subset S[t, t^{-1}] .
\]
Here, we think that $S[t, t^{-1}]$ is a $\bZ^{2}$-graded ring with
$\deg(x) = (0,a)$, $\deg(y) = (0,b)$, $\deg(z) = (0,c)$, $\deg(t) = (1,0)$.
Then $A$ is a $\bZ^{2}$-graded subring of $S[t, t^{-1}]$.

In this situation, $A$ is Noetherian if and only if it satisfies Huneke's criterion~\cite{Hu}, that is, there exist positive integers $r$ and $s$, 
and elements $f \in P^{(r)}$ and $g \in P^{(s)}$ such that
\begin{equation}\label{huneke}
\ell (S/(f,g,x)) = rs a .
\end{equation}
One can prove that, once $A$ is Noetherian,
we can choose homogeneous elements $f$ and $g$ satisfying (\ref{huneke})
as in Cutkosky~\cite{Cu}.
Assume that $A$ is Noetherian and homogeneous elements $f$ 
and $g$ satisfy Huneke's criterion.
Put $\deg(f) = (0,d_{1})$ and $\deg(g) = (0,d_{2})$.
In this case, we can prove that
\[
B = k[t^{-1}, x, ft^{r}, gt^{s}] \subset A
\]
is a finite map.
Then, by Proposition~\ref{reduce} (4), ray ideal cones of $A$ are the same as those of $B$.
Remark that $t^{-1}$, $x$, $ft^{r}$, $gt^{s}$ are algebraically independent
over $k$.
Degrees of them are $(-1,0)$, $(0,a)$, $(r,d_{1})$, $(s,d_{2})$, respectively.
Here assume that $d_{1}/r < d_{2}/s$ is satisfied~\footnote{
In many cases,  $d_{1}/r$ does not equal $d_{2}/s$
when there exist homogeneous elements $f$, $g$ satisfying (\ref{huneke}).
In the case $(a,b,c) = (1,1,1)$,
$d_{1}/r = d_{2}/s = 1$ holds.
The authors do not know any other examples satisfying  $d_{1}/r = d_{2}/s$.
}.
Then, by Example~\ref{2.2}, the non-zero ray ideals of $B$ are 
\begin{align*}
& B, \\
& (t^{-1})B, \ \ \ (t^{-1}, x)B \cap (x, gt^{s}, ft^{r})B, \ \ \ (t^{-1}, x, gt^{s})B \cap ( gt^{s}, ft^{r}), \ \ \
(ft^{r})B, \\
& (t^{-1})B \cap (x, gt^{s}, ft^{r})B, \ \ \ (t^{-1}, x)B \cap (gt^{s}, ft^{r})B, \ \ \ (t^{-1},x,gt^s)B \cap (ft^{r})B .
\end{align*}
The corresponding maximal ray ideal cones are
\begin{align*}
& \{ \bmo \}, \\
& \bRo   (-1,0), \ \ \ \bRo   (0,1), \ \ \ \bRo  (s,d_{2}), \ \ \
\bRo  (r,d_{1}), \\
& \bRo   (-1,0)+\bRo   (0,1), \ \ \ \bRo   (0,1)+\bRo  (s,d_{2}), \ \ \ \bRo  (s,d_{2})+\bRo  (r,d_{1}) .
\end{align*}
One can see that Theorem~\ref{eti} is satisfied in this case.

Next, we assume that $(a,b,c)=(25,29,72)$.
Then the symbolic Rees ring $A$ is not Notherian by Goto-Nishida-Watanabe~\cite{GNW}.
There exist a homogeneous element $ft \in A_{(1, 216)}$ and an ideal $J$ of $A$ such that,
for positive rational numbers $\alpha$ and $\beta$, we have
\[
J_{(\alpha, \beta)}(A) = \left\{
\begin{array}{lll}
J & & \mbox{if $\beta/\alpha > 25{\cdot} 29/3$,} \\
(ft)A \cap J & & \mbox{if $25{\cdot} 29/3 \ge \beta/\alpha > 216$,} \\
(ft)A & & \mbox{if $\beta/\alpha = 216$,} \\
0 & & \mbox{if $216 > \beta/\alpha $.}
\end{array}
\right.
\]
Here remark that
\[
J \supsetneq (ft)A \cap J \subsetneq (ft)A .
\]
Thus Theorem~\ref{eti} is not true if we remove the assumption that $A$ is Noetherian.
Since the height of $J$ is bigger than $1$,
this example satisfies the condition (I) in Theorem~\ref{Demazure}.
\end{rm}
\end{Example}

\section{When is $(\bRo)^{n}$ a chamber?}

The purpose of this section is to prove the following theorem:

\begin{Theorem}\label{Arai}
Let $A$ be a Noetherian $\bZ^{n}$-graded domain such that $A_{\bmo}$ is a field.
Assume that $C(A) = (\bRo)^{n}$.
Let $\bme_{i}$ be the $i$-th unit vector in $\bR^{n}$.
Let $B$ be the subring of $A$ generated by
\[
\bigcup_{i=1}^{n} \left[ \bigcup_{m > 0}A_{m \bme_{i}}\right]
\]
over $A_{\bmo}$.

Then the following conditions are equivalent:
\begin{enumerate}
\item
The inclusion $B \rightarrow A$ is a finite extension.
\item
Any face of $(\bRo)^{n}$ is a ray ideal cone of $A$.
\item
$C(A)$ itself is a ray ideal cone of $A$.
\end{enumerate}
\end{Theorem}

\proof
First, we shall prove $(1) \Rightarrow (2)$.
By Proposition~\ref{reduce} (4), we may assume that $B = A$.
Restricting $\bZ^{n}$ to a subgroup generated by a subset of $\{ \bme_1, \ldots, \bme_n\}$
 (cf.\ Proposition~\ref{reduce} (2)),
we have only to show that $C(A)$ itself is a ray ideal cone of $A$.
Take
\[
\bma = a_{1} \bme_{1} + a_{2} \bme_{2} + \cdots + a_{n} \bme_{n} \in {\rm int}(C(A)) \cap \bZ^{n} ,
\]
where $a_{i}$'s are positive integers.
Take $x_{i} \in A_{b_{i}\bme_{i}}$ for $i = 1, 2, \ldots, n$, 
where $b_{i}$'s are positive integers.
Put $b = b_{1}b_{2} \cdots b_{n}$.
Then,
\[
\deg(x_{1}^{(ba_{1}/b_{1})}x_{2}^{(ba_{2}/b_{2})} \cdots x_{n}^{(ba_{n}/b_{n})})
= b \bma .
\]
Therefore we know $J_{\bma}(A) \supset J_{\bme_{1}}(A)J_{\bme_{2}}(A)\cdots J_{\bme_{n}}(A)$.
On the other hand, since $A$ is generated by homogeneous elements with degree
on coordinate axes,
we have $I_{\bma}(A) \subset I_{\bme_{1}}(A)I_{\bme_{2}}(A)\cdots I_{\bme_{n}}(A)$
immediately.
Thus we obtain 
\[
J_{\bma}(A) = \sqrt{J_{\bme_{1}}(A)J_{\bme_{2}}(A)\cdots J_{\bme_{n}}(A)} .
\]
Therefore $C(A)$ itself is a ray ideal cone of $A$.

The implication $(2) \Rightarrow (3)$ is trivial.

In the rest of this section, we shall prove $(3) \Rightarrow (1)$.
There exists a positive integer $c$ such that $A_{c\bme_{i}} \neq 0$ for $i = 1, 2, ,\ldots, n$.
Replacing $A$ by $A_{c\bZ^{n}}$, we may assume $A_{\bme_{i}} \neq 0$ for $i = 1, 2, ,\ldots, n$.
Let $\tilde{A}$ be the normalization of $A$.
Since $A$ is finitely generated over the field $A_{\bmo}$,
$\tilde{A}$ is also a finitely generated $\bZ^{n}$-graded ring
over the field $\tilde{A}_{\bmo}$.
Then $C(\tilde{A}) = C(A) = (\bRo)^{n}$, and
$(\bRo)^{n}$ is a chamber of $\tilde{A}$ by Proposition~\ref{reduce} (4).
Let $B'$ be the subring of $\tilde{A}$ generated by
\[
\bigcup_{i=1}^{n} \left[ \bigcup_{m > 0}\tilde{A}_{m \bme_{i}}\right]
\]
over $\tilde{A}_{\bmo}$.
Consider the following diagram:
\[
\begin{array}{ccc}
B & \longrightarrow & A \\
\downarrow & & \downarrow \\
B' & \longrightarrow & \tilde{A}
\end{array}
\]
It is easy to see that $B \rightarrow B'$ is a finite extension.
If $B' \rightarrow \tilde{A}$ is a finite extension,
then so is $B \rightarrow A$.
Replacing $A$ with $\tilde{A}$, we may assume that $A$ is a normal domain.

We denote the field of fractions by $Q( \ )$.

Let $\overline{B}$ be the integral closure of $B$ in $A$.
Since $Q(A)$ is finitely generated over $Q(B)$ as a field,
$\overline{B}$ is finite over $A$.
We shall prove $A=\overline{B}$.

\vspace{2mm}

Let us prove the following claim:

\begin{Claim}\label{claimA1}
The ray ideal $J_{C(A)}$ is $\sqrt{J_{\bme_1}(A) J_{\bme_2}(A) \cdots J_{\bme_n}(A)}$.
\end{Claim}

It is easy to see that  $J_{C(A)}$ contains $\sqrt{J_{\bme_1}(A) J_{\bme_2}(A) \cdots J_{\bme_n}(A)}$ as in the proof of the implication $(1) \Rightarrow (2)$.
We shall prove the opposite inclusion.
Let $P$ be a minimal prime ideal of $\sqrt{J_{\bme_1}(A) J_{\bme_2}(A) \cdots J_{\bme_n}(A)}$.
It is well known that $P$ is $\bZ^n$-graded.
We want to show that $P$ contains $J_{C(A)}$.
There exists $i$ such that $P$ contains $J_{\bme_i}(A)$.
Then $C(A/P)$ does not contain $\bme_i$.
Since $A$ is Noetherian, both $C(A)$ and $C(A/P)$ are rational polyhedral cones.
Since $C(A/P) \subsetneq C(A)$, there exists $\bmb \in {\rm int}(C(A)) \cap \bQ^{n}$ such that 
$\bmb \not\in C(A/P)$.
Therefore $P$ contains $J_{\bmb}(A) = J_{C(A)}$.
We have completed the proof of Claim~\ref{claimA1}.

\vspace{2mm}

Next, we shall prove the following claim:

\begin{Claim}\label{claimA2}
For $\bma \in \bN^{n}$, $A_{\bma}$ is contained in $\overline{B}$.
\end{Claim}

There exists $b \in \bZ$ such that $J_{\bme_i}(A) = \sqrt{A_{b\bme_i}A}$ for $i = 1, 2, \ldots, n$, where
$A_{b\bme_i}A$ denotes the ideal of $A$ generated by $A_{b\bme_i}$.

We put $\bma = a_1\bme_1 + a_2\bme_2 + \cdots + a_n\bme_n$, where $a_1, a_2, \ldots, a_n \in \bN$.
Then, by Claim~\ref{claimA1}, we have
\begin{eqnarray*}
J_\bma(A) & = & \sqrt{J_{\bme_1}(A) J_{\bme_2}(A) \cdots J_{\bme_n}(A)} \\
& = & \sqrt{A_{b\bme_1}A_{b\bme_2}\cdots A_{b\bme_n}A} \\
& = & \sqrt{A_{ba_1\bme_1}A_{ba_2\bme_2}\cdots A_{ba_n\bme_n}A} ,
\end{eqnarray*}
where $A_{ba_1\bme_1}A_{ba_2\bme_2}\cdots A_{ba_n\bme_n}A$ denotes the ideal of $A$ generated by 
\[
A_{ba_1\bme_1}A_{ba_2\bme_2}\cdots A_{ba_n\bme_n} =
\{ x_{1}x_{2}\cdots x_{n} \mid x_{1} \in A_{ba_1\bme_1}, \, x_{2} \in A_{ba_2\bme_2}, \, \ldots,\,  x_{n} \in A_{ba_n\bme_n} \} .
\]
Then there exists $m \in \bN$ such that
\begin{equation}\label{siki1}
J_\bma(A)^{m} \subset A_{ba_1\bme_1}A_{ba_2\bme_2}\cdots A_{ba_n\bme_n}A .
\end{equation}
We put
\[
T_{\ell} = A_{\ell b \bma} \ \ \mbox{and} \ \ T = \bigoplus_{\ell \in \bNo} T_{\ell} .
\]
Since $T$ is finitely generated over $T_{0} = A_{\bmo}$ as a ring,
we have
\begin{equation}\label{siki2}
T_{\ell} \subset J_\bma(A)^{m}
\end{equation}
for $\ell \gg 0$.
Here, recall that 
\[
A_{ba_1\bme_1}A_{ba_2\bme_2}\cdots A_{ba_n\bme_n} \subset T_{1} .
\]
Then, by (\ref{siki1}) and (\ref{siki2}),
\[
T_{\ell} = A_{ba_1\bme_1}A_{ba_2\bme_2}\cdots A_{ba_n\bme_n} T_{\ell -1}
\]
for $\ell \gg 0$.
Then we know that the inclusion
\[
A_{\bmo}[A_{ba_1\bme_1}A_{ba_2\bme_2}\cdots A_{ba_n\bme_n}]
\longrightarrow
T
\]
is a finite morphism.
Since $B$ contains $A_{\bmo}[A_{ba_1\bme_1}A_{ba_2\bme_2}\cdots A_{ba_n\bme_n}]$,
$\overline{B}$ contains $T$.
Since elements in $A_{\bma}$ are integral over $T$, $A_{\bma}$ is contained in $\overline{B}$.
We have completed the proof of Claim~\ref{claimA2}.

\vspace{2mm}

Now, we start to prove $\overline{B} = A$.
It is enough to prove $\overline{B}_{P} \supset A$
for a prime ideal $P$ of $\overline{B}$ of height one
since $\overline{B}$ is a normal domain.
By Claim~\ref{claimA2}, $Q(\overline{B}) = Q(A)$.
Let $v_{P}$ be the valuation associated to the
discrete valuation ring $\overline{B}_{P}$.

First, assume that, for each $i$, there exists $x_{i} \in A_{a_{i}\bme_{i}}$ for some $a_{i} \in \bN$
such that $v_{P}(x_{i}) = 0$.
For any non-zero homogeneous element $x$ in $A$, 
\[
\deg(x x_{1} x_{2} \cdots x_{n}) \in \bN^{n} .
\]
Then, by Claim~\ref{claimA2}, we have
\[
x x_{1} x_{2} \cdots x_{n} \in \overline{B}  .
\]
Then we have $v_{P}(x) \ge 0$ immediately.
Thus $x$ is in $\overline{B}_{P}$ in this case.

Next, assume that there exists $i$ such that $P$ contains
\[
\bigcup_{m \in \bN}A_{m \bme_i} .
\]
For the simplicity, suppose that the above $i$ is $1$.
Then $P \cap B$ contains all the homogeneous element such that the first component of the degree is positive.
Since $\overline{B}$ is integral over $B$, $P$ contains all the homogeneous element in $\overline{B}$
such that the first component of the degree is positive.
Let $\overline{B}_{(0)}$ be the homogeneous localization of $\overline{B}$, that is
\[
\overline{B}_{(0)} =
\left\{ \left. z \in Q(\overline{B}) \ \right| \ \mbox{$yz \in \overline{B}$ for some non-zero homogeneous element
$y \in \overline{B}$ } \right\} .
\]
It is still a $\bZ^{n}$-graded ring.
Then, by Claim~\ref{claimA2}, it is easy to see $A \subset \overline{B}_{(0)}$.
Since $A_{\bme_{i}} \neq 0$ for $i = 1, 2, \ldots, n$,
we may assume that $\overline{B}_{(0)} = K[t_{1}^{\pm 1}, t_{2}^{\pm 1}, \ldots, t_{n}^{\pm 1}]$,
where $\deg(t_{i}) = \bme_{i}$ and $K$ is a field such that $(\overline{B}_{(0)})_{\bmo} = K$.
Then $P = t_{1}K[t_{1}, t_{2}, \ldots, t_{n}] \cap \overline{B}$ 
since the height of $P$ is $1$.
Hence
$\overline{B}_{P} = K[t_{1}, t_{2}, \ldots, t_{n}]_{(t_{1})}$
because $\overline{B}_{P}$ is a discrete valuation ring.
Thus we have
\[
A \subset K[t_{1}, t_{2}, \ldots, t_{n}] \subset K[t_{1}, t_{2}, \ldots, t_{n}]_{(t_{1})} =
\overline{B}_{P} .
\]

We have completed the proof of
Theorem~\ref{Arai}.
\qed

\section{Chamber decomposition}

The aim of this section is to prove the following theorem, which also follows from
Theorem~2.11 in \cite{BH}. 

\begin{Theorem}\label{ChamberDecomp}
Let $A$ be a Noetherian $\bZ^{n}$-graded domain such that $A_{\bmo}$ is a field.
Then there exists finitely many chambers $\sigma_{1}$, $\sigma_{2}$, \ldots, $\sigma_{m}$
such that 
\[
C(A) = \bigcup_{i=1}^{m} \sigma_{i} .
\]
\end{Theorem}

\proof
Replacing $\bZ^{n}$ by the subgroup generated by $\{ \bma \in \bZ^{n} \mid A_{\bma} \neq 0 \}$, we may assume that $\dim C(A) = n$.

Let $x_{1}$, $x_{2}$, \ldots, $x_{s}$ be a set of non-zero homogeneous elements such that
$A = A_{\bmo}[x_{1}, x_{2}, \ldots, x_{s}]$. 
Put $\bma_{i} = \deg(x_{i})$ for each $i$.
Then,
\[
C(A) = \sum_{i=1}^s\bRo   \bma_{i} .
\]
Here, put
\[
\{ H_{1}, H_{2}, \ldots, H_{\ell} \} \\
= 
\left\{ H \ \left| 
\begin{array}{l}
\mbox{$H$ is an $(n-1)$-dimensional linear subspace of $\bR^{n}$} \\ 
\mbox{spanned by a subset of $\bma_{1}$, $\bma_{2}$, \ldots, $\bma_{s}$} 
\end{array}
\right. \right\} .
\]
Let $f_{i}: \bR^{n} \rightarrow \bR$ be a $\bR$-linear map such that
\[
H_{i} = \{ \bma \in \bR^{n} \mid f_{i}(\bma) = 0 \} 
\]
for each $i$.
For $\xi_{1}, \xi_{2}, \ldots, \xi_{\ell} \in \{ -1, 1 \}$, we put
\[
C(\xi_{1}, \xi_{2}, \ldots, \xi_{\ell})
= \{ \bma \in \bR^{n} \mid \mbox{$\xi_{i} f_{i}(\bma) > 0$ for $i = 1, 2, \ldots, \ell$} 
\} .
\]

\vspace{2mm}

Here, we prove the following claim:

\vspace{2mm}

\begin{Claim}~\label{claimE1}
Let $\overline{C(\xi_{1}, \xi_{2}, \ldots, \xi_{\ell})}$
be the closure of $C(\xi_{1}, \xi_{2}, \ldots, \xi_{\ell})$
in the classical topology of $\bR^{n}$.
Then
\begin{equation}\label{5.1}
C(A) = \bigcup_{C(A) \cap C(\xi_{1}, \xi_{2}, \ldots, \xi_{\ell}) \neq \emptyset}
\overline{C(\xi_{1}, \xi_{2}, \ldots, \xi_{\ell})} .
\end{equation}
\end{Claim}

First, we prove that the right-hand-side is contained in the left one.
Assume that $C(A) \cap C(\xi_{1}, \xi_{2}, \ldots, \xi_{\ell}) \neq \emptyset$.
Take $\bmx \in C(A) \cap C(\xi_{1}, \xi_{2}, \ldots, \xi_{\ell})$.
Since 
\[
\bmx \in C(A) = \sum_{i}\bRo   \bma_{i} ,
\]
there exist linearly independent vectors $\bma_{i_{1}}$, $\bma_{i_{2}}$, \ldots, $\bma_{i_{n}}$ such that 
\[
\bmx \in \bRo   \bma_{i_{1}}+ \bRo   \bma_{i_{2}}+ \cdots + \bRo   \bma_{i_{n}}
\]
by Carath\'eodory's theorem.
Since $\bmx \in C(\xi_{1}, \xi_{2}, \ldots, \xi_{\ell})$,
$\bmx$ is not contained in hypersurface spanned by a subset of $\bma_{i_{1}}$, $\bma_{i_{2}}$, \ldots, $\bma_{i_{n}}$.
Hence we know 
\begin{equation}\label{x}
\bmx \in \bRp   \bma_{i_{1}}+ \bRp   \bma_{i_{2}}+ \cdots + \bRp   \bma_{i_{n}} .
\end{equation}
Take any $\bmy \in C(\xi_{1}, \xi_{2}, \ldots, \xi_{\ell})$.
Since $\bma_{i_{1}}$, $\bma_{i_{2}}$, \ldots, $\bma_{i_{n}}$ span $\bR^{n}$,
we can write
\begin{equation}\label{y}
\bmy = r_{1} \bma_{i_{1}} + r_{2} \bma_{i_{2}} + \cdots + r_{n}\bma_{i_{n}} ,
\end{equation}
where $r_{1}, r_{2}, \ldots, r_{n} \in \bR$.
Assume that $r_{j} \le 0$ for some $j$.
Let $H_{t}$ be the hypersurface spanned by 
$\bma_{i_{1}}$, \ldots, $\bma_{i_{j-1}}$, $\bma_{i_{j+1}}$, \ldots, $\bma_{i_{n}}$.
By (\ref{x}), we obtain $\xi_{t}f_{t}(\bma_{i_{j}}) > 0$.
On the other hand, by (\ref{y}), 
$0 < \xi_{t}f_{t}(\bmy) = r_{j}\xi_{t}f_{t}(\bma_{i_{j}})$.
It is a contradiction.
Hence we obtain
\begin{equation}\label{z}
C(\xi_{1}, \xi_{2}, \ldots, \xi_{\ell})
\subset \bRp   \bma_{i_{1}}+ \bRp   \bma_{i_{2}}+ \cdots + \bRp   \bma_{i_{n}} \subset C(A) .
\end{equation}
Since $C(A)$ is a closed subset, the closure $\overline{C(\xi_{1}, \xi_{2}, \ldots, \xi_{\ell})}$
is contain in $C(A)$.

Next, we shall show that the right-hand-side contains the left one in (\ref{5.1}).
It is enough to show that any point $\bmz$ in the interior of $C(A)$ is 
contained in the right-hand-side.
Let $U_{\epsilon}(\bmz)$ be an open ball of radius $\epsilon$ with center $\bmz$.
If $\epsilon$ is small enough, then $U_{\epsilon}(\bmz)$ is contained in $C(A)$.
Since
\[
U_{\epsilon}(\bmz) \not\subset H_{1} \cup H_{2} \cup \cdots \cup H_{\ell} ,
\]
$U_{\epsilon}(\bmz)$ intersects the right-hand-side in (\ref{5.1}).
Since the right-hand-side is a closed set, $\bmz$ is contained in the right-hand-side.

We have completed the proof of Claim~\ref{claimE1}.

\vspace{2mm}

In order to prove Theorem~\ref{ChamberDecomp}, it is enough to show the following claim:

\vspace{2mm}

\begin{Claim}~\label{claimE2}
Suppose $C(A) \cap C(\xi_{1}, \xi_{2}, \ldots, \xi_{\ell}) \neq \emptyset$.
Then $\overline{C(\xi_{1}, \xi_{2}, \ldots, \xi_{\ell})}$ is a chamber of $A$.
\end{Claim}

Suppose that $C(A) \cap C(\xi_{1}, \xi_{2}, \ldots, \xi_{\ell}) \neq \emptyset$.
Then it is easy to see that
\[
\overline{C(\xi_{1}, \xi_{2}, \ldots, \xi_{\ell})}  = \{ \bmx \in \bR^{n} \mid \forall i, \ \ \xi_{i}f_{i}(\bmx) \ge 0  \} 
\]
is a rational polyhedral cone.
The interior of $\overline{C(\xi_{1}, \xi_{2}, \ldots, \xi_{\ell})}$ is
$C(\xi_{1}, \xi_{2}, \ldots, \xi_{\ell})$.
Suppose $\deg(x_{1}^{t_{1}} x_{2}^{t_{2}} \cdots x_{s}^{t_{s}}) \in C(\xi_{1}, \xi_{2}, \ldots, \xi_{\ell})$.
We shall prove that $x_{1}^{t_{1}} x_{2}^{t_{2}} \cdots x_{s}^{t_{s}} \in J_{\bma}(A)$
for any $\bma \in C(\xi_{1}, \xi_{2}, \ldots, \xi_{\ell}) \cap \bQ^{n}$.
Since
\[
\deg(x_{1}^{t_{1}} x_{2}^{t_{2}} \cdots x_{s}^{t_{s}}) \in \sum_{t_i > 0} \bRo   \bma_i ,
\]
there exists $i_{1}$, $i_{2}$, \ldots, $i_{n}$ such that
\begin{itemize}
\item
$\bma_{i_{1}}$, $\bma_{i_{2}}$, \ldots, $\bma_{i_{n}}$ are linearly independent, 
\item
$t_{i_{1}} > 0$, $t_{i_{2}} > 0$, \ldots, $t_{i_{n}} > 0$, and
\item
$\deg(x_{1}^{t_{1}} x_{2}^{t_{2}} \cdots x_{s}^{t_{s}})
\in \bRp   \bma_{i_{1}}+ \bRp   \bma_{i_{2}}+ \cdots + \bRp   \bma_{i_{n}}$
\end{itemize}
by the same argument just before (\ref{x}).
Then, by the argument that (\ref{x}) implies (\ref{z}),
we obtain 
\[
\bma \in C(\xi_{1}, \xi_{2}, \ldots, \xi_{\ell}) \subset
\bRp   \bma_{i_{1}}+ \bRp   \bma_{i_{2}}+ \cdots + \bRp   \bma_{i_{n}} .
\]
Then it is easy to see that a power of $x_{1}^{t_{1}} x_{2}^{t_{2}} \cdots x_{s}^{t_{s}}$ is contained
in $I_{\bma}(A)$.

We have completed the proofs of Claim~\ref{claimE2} and Theorem~\ref{ChamberDecomp}.
\qed

\section{The fan structure of ray ideal cones}

In this section, we shall prove that, for a Noetherian $\bZ^{n}$-graded domain, 
there exists the unique 
maximal ray ideal cone $\sigma_J$ for each non-zero ray ideal $J$.
Furthermore, we shall prove that the set of maximal ray ideal cones forms a fan
in Theorem~\ref{FanStructure}, which follows from Theorem~2.11 in \cite{BH}.
If $A$ is not Notherian, Theorem~\ref{FanStructure} is not true.
In fact, the set of maximal ray ideal cones of the non-Noetherian 
symbolic Rees ring in Example~\ref{SpacceMonomial} do not form a fan.

\begin{Lemma}\label{face}
Let $A$ be a Noetherian $\bZ^n$-graded ring and $\sigma$ be a ray ideal cone of $A$.
Take $\bmb \in \sigma \cap \bQ^n$.
Then $J_{\bmb}(A)$ contains the ray ideal $J_{\sigma}$ of the ray ideal cone $\sigma$.
\end{Lemma}

\proof
Let $P$ be a homogeneous prime ideal containing $J_{\bmb}(A)$.
Remark that $C(A/P)$ is a closed subset of $\bR^{n}$
since $A$ is Noetherian.
If $C(A/P)$ contains ${\rm rel.int}(\sigma) \cap \bQ^n$, then
$C(A/P)$ contains $\sigma$.
However, $C(A/P)$ does not contain $\bmb$.
Therefore there exists  $\bma \in {\rm rel.int}(\sigma) \cap \bQ^n$
such that $P \supset J_{\bma}(A) = J_{\sigma}$.
Thus $J_{\bmb}(A)$ contains $J_{\sigma}$.
\qed

\begin{Theorem}\label{FanStructure}
Let $A$ be a Noetherian $\bZ^{n}$-graded domain such that $A_{\bmo}$ is a field.
\begin{enumerate}
\item
There exists only finitely many non-zero ray ideals.
\item
For each non-zero ray ideal $J$, there exists the unique maximal ray ideal cone $\sigma_{J}$ of the ray ideal $J$.
Furthermore, $\sigma_{J}$ is a rational polyhedral cone.
For any $\bma \in \bQ^n \setminus \sigma_J$,
$J_{\bma}(A)$ does not coincide with $J$.
\item
Suppose that a rational polyhedral cone $\sigma$ is a ray ideal cone.
Then any face of $\sigma$ is also a ray ideal cone.
\item
Let $J_{1}$ and $J_{2}$ be non-zero ray ideals of $A$.
Then $J_{1}$ contains $J_{2}$ if and only if $\sigma_{J_{1}}$ is a face of $\sigma_{J_{2}}$.
\item
Assume that $A_{\bmo} \setminus \{ 0 \}$ coincides with the set of units in $A$.
Then
\[
\{ \sigma_{J} \mid \mbox{$J$ is a non-zero ray ideal} \}
\]
forms a fan.
\end{enumerate}
\end{Theorem}

\proof
First, we shall prove (1).
It is enough to show that there exist finitely many ray ideal cones $\rho_{1}$, $\rho_{2}$, \ldots, $\rho_{\ell}$ such that
\begin{equation}\label{decomp}
C(A) = \bigcup_{i = 1}^{\ell}{\rm rel.int}(\rho_{i}) ,
\end{equation}
where each $\rho_{i}$ is a rational polyhedral cone.
(If it is satisfied, $A$ has only finitely many non-zero ray ideals
$J_{\rho_{1}}$, $J_{\rho_{2}}$, \ldots, $J_{\rho_{\ell}}$.)
We shall prove it by induction on the dimension of $C(A)$.
If $\dim C(A) = 0$, it is trivial.
Assume that $\dim C(A) > 0$.
By Theorem~\ref{ChamberDecomp}, we have chambers
$\sigma_1$, $\sigma_2$, \ldots, $\sigma_m$ satisfying
\[
C(A) = \bigcup_{i = 1}^{m} \sigma_{i} .
\]
It is sufficient to show (\ref{decomp}) for $A_{\sigma_{1}}$, $A_{\sigma_{2}}$, \ldots, $A_{\sigma_{m}}$ since any ray ideal cone of $A_{\sigma_i}$ is a ray ideal cone of $A$
by Proposition~\ref{reduce} (1).
Replacing $A$ by $A_{\sigma_{i}}$, we assume that $C(A)$ itself is a chamber of $A$ and
$C(A)$ is a rational polyhedral cone.
Then the boundary of $C(A)$ is a union of finitely many faces of $C(A)$.
Let $\tau$ be a face of $C(A)$.
Since $\dim \tau < \dim C(A)$,
(\ref{decomp}) holds for $A_{\tau}$.
Therefore (\ref{decomp}) is satisfied
since any ray ideal cone of $A_\tau$ is a ray ideal cone of $A$
by Proposition~\ref{reduce} (1).

Next, we shall prove (2).
Let $J$ be a non-zero ray ideal of $A$.
Consider a decomposition of $C(A)$ as in (\ref{decomp}).
Suppose that the ray ideal of $\rho_{i}$ is $J$ if and only if $i = 1, 2, \ldots, s$.
Then it is easy to see
\begin{equation}\label{3siki}
\{ \bma \in \bQ^{n} \mid J_{\bma}(A) = J \} 
\subset \bigcup_{i = 1}^{s}{\rm rel.int}(\rho_{i})
\subset \sum_{\bma \in \bQ^{n}, \ J_{\bma}(A) = J} \bRo   \bma .
\end{equation}
Take an element
\[
\bmx = u_{1} \bma_{1} + u_{2} \bma_{2} + \cdots + u_{t} \bma_{t}
\]
in the right end of (\ref{3siki}), where $u_{i} \in \bRp$ and $J_{\bma_{i}}(A) = J$ for $i = 1, 2, \ldots, t$.
Let $\{ u_{ij} \}_{j}$ be a sequence of positive rational numbers which converges to $u_{i}$.
We put
\[
\bmx_{j} = u_{1j} \bma_{1} + u_{2j} \bma_{2} + \cdots + u_{tj} \bma_{t}
\]
for each $j$.
By Theorem~\ref{eti}, $J_{\bmx_{j}}(A) =J$ for all $j$.
Since $\{ \bmx_{j} \}_{j}$ converges to $\bmx$,
$\bmx$ is contained in the closure of the left end in (\ref{3siki}).
Taking the closure of each set in (\ref{3siki}), we obtain
\begin{equation}\label{closure}
\overline{ \{ \bma \in \bQ^{n} \mid J_{\bma}(A) = J \} }
= \bigcup_{i = 1}^{s} \rho_{i}
= \overline{ \sum_{\bma \in \bQ^{n}, \ J_{\bma}(A) = J} \bRo \bma } .
\end{equation}
The right end of (\ref{closure}) is a cone in $\bR^{n}$.
We denote this cone by $\sigma_{J}$.
Since $\sigma_{J}$ is a cone spanned by $\rho_{1}$, $\rho_{2}$, \ldots, $\rho_{s}$,
it is a rational polyhedral cone.
It is easy to see that $\sigma_{J}$ is the unique maximal ray ideal cone with ray ideal $J$.

Next, we shall prove (3).
If $\dim \sigma \le 2$, then it is easy.
Assume that $\dim \sigma \ge 3$.
Assume that $\tau$ is a proper face of $\sigma$ such that
there exist $\bma, \bmb \in {\rm rel.int}(\tau) \cap \bQ^{n}$ satisfying
$J_{\bma}(A) \neq J_{\bmb}(A)$.
We can take $\bma', \bmb' \in {\rm rel.int}(\tau) \cap \bQ^{n}$ such that
$\bma, \bmb \in {\rm rel.int}(\bRo   \bma' + \bRo   \bmb')$.
Take $\bmc \in {\rm rel.int}(\sigma) \cap \bQ^{n}$.
Put
\[
T = \bRo   \bma' + \bRo   \bmb' + \bRo   \bmc .
\]
By definition, $\dim T = 3$.
The relative interior of $T$ is contained in that of $\sigma$.
By Proposition~\ref{reduce} (3), $T$ itself is a ray ideal cone of $A_{T}$.
Then, by Theorem~\ref{Arai}, any face of $T$ is a ray ideal cone of $A_T$.
Since $\bma, \bmb \in {\rm rel.int}(\bRo   \bma' + \bRo   \bmb')$, we know $J_{\bma}(A_{T}) = J_{\bmb}(A_{T})$.
By Proposition~\ref{reduce} (1), we have $J_{\bma}(A) = J_{\bmb}(A)$.
It is a contradiction.

Next, we shall prove (4).
If $\sigma_{J_{1}}$ is a face of $\sigma_{J_{2}}$,
then $J_{1}$ contains $J_{2}$ by Lemma~\ref{face}.

Conversely assume that $J_{1} \supsetneq J_{2}$.
Let $\bma \in {\rm rel.int}(\sigma_{J_{2}}) \cap \bQ^{n}$
and $\bmb \in {\rm rel.int}(\sigma_{J_{1}}) \cap \bQ^{n}$.
Let $\bmc$ be a rational point on the line segment
between $\bma$ and $\bmb$ such that $\bmc \neq \bmb$.
Then $J_{\bmc}(A) = J_{\bma}(A)$ by Theorem~\ref{eti}.
Thus one can see
\[
{\rm rel.int}(\sigma_{J_{1}}) \subset \overline{{\rm rel.int}(\sigma_{J_{2}})} .
\]
Hence 
\[
\sigma_{J_{1}} \subset \sigma_{J_{2}} .
\]
Since $\sigma_{J_{1}}$ is contained in the boundary of $\sigma_{J_{2}}$,
there exists the minimal face $\tau$ of $\sigma_{J_{2}}$ containing $\sigma_{J_{1}}$.
By the minimality of $\tau$, $\sigma_{J_{1}}$ is not contained in the boundary of $\tau$.
Here, remark that
\[
\sigma_{J_1} = \overline{ \{ \bma \in \bQ^{n} \mid J_{\bma}(A) = J_1 \} } 
\]
by (\ref{closure}).
Therefore there exists $\bmd \in {\rm rel.int}(\tau) \cap \bQ^n$
such that $J_{\bmd}(A) = J_1$.
By (3), $\tau$ is a ray ideal cone of the ray ideal $J_{1}$.
Since $\sigma_{J_{1}}$ is the unique maximal ray ideal cone with ray ideal $J_{1}$,
$\sigma_{J_{1}}$ coincides with $\tau$.

Next, we shall prove (5).
Put
\[
\Delta(A) = \{ \sigma_{J} \mid \mbox{$J$ is a non-zero ray ideal} \} .
\]
Since $A_{\bmo} \setminus \{ 0 \}$ coincides with the set of units in $A$,
each $\sigma_J$ in $\Delta(A)$ is a strongly convex rational polyhedral cone.
In order to show that $\Delta(A)$ forms a fan, it is enough to show
\begin{itemize}
\item[(a)]
any face of a cone in $\Delta(A)$ is in $\Delta(A)$, and
\item[(b)]
for $\sigma_1, \sigma_2 \in \Delta(A)$, $\sigma_1 \cap \sigma_2$
is a face of both $\sigma_1$ and $\sigma_2$.
\end{itemize}

First, we shall prove (a).
Take $\sigma_J \in \Delta(A)$.
Let $\tau$ be a face of $\sigma_J$.
By (3), $\tau$ is a ray ideal cone of $A$.
Let $J'$ be the ray ideal of the ray ideal cone $\tau$.
By Lemma~\ref{face}, $J'$ contains $J$.
Then, by (4), $\sigma_{J'}$ is a face of $\sigma_J$.
Since $\sigma_{J'}$ contains $\tau$,
$\tau$ is a face of $\sigma_{J'}$.
We want to show $\tau = \sigma_{J'}$.
Assume that $\tau \subsetneq \sigma_{J'}$.
Take $\bma \in {\rm rel.int}(\sigma_{J'}) \cap \bQ^n$ and
$\bmb \in {\rm rel.int}(\tau) \cap \bQ^n$.
Then we have $J_{\bma}(A) = J_{\bmb}(A) \neq 0$.
By Corollary~\ref{open}, there exists $\bmc \in \bQ^n$
such that $J_{\bmc}(A) = J_{\bmb}(A)$ and $\bmc \not\in \sigma_{J'}$.
It contradicts the maximality of $\sigma_{J'}$ in (2).
Hence we have $\tau = \sigma_{J'}$.
Therefore any face of $\sigma_J$ is contained in $\Delta(A)$.

Next, we shall prove (b).
Take $\sigma_{J_1}, \sigma_{J_2} \in \Delta(A)$.
We shall prove that $\sigma_{J_1} \cap \sigma_{J_2}$
is a face of both $\sigma_{J_1}$ and $\sigma_{J_2}$.
If $J_1 \supset J_2$ or $J_1 \subset J_2$,
then the assertion follows from (4).
Assume $J_1 \not\supset J_2$ and $J_1 \not\subset J_2$.
Since both $\sigma_{J_1}$ and $\sigma_{J_2}$ are rational 
polyhedral cones, so is $\sigma_{J_1} \cap \sigma_{J_2}$.
By Lemma~\ref{face}, $J_{\bma}(A)$ contains both $J_1$ and $J_2$ for any $\bma \in \sigma_{J_1} \cap \sigma_{J_2} \cap \bQ^n$.
Then  $\sigma_{J_1} \cap \sigma_{J_2}$ is contained in 
the boundary of $\sigma_{J_i}$ for $i = 1, 2$ by (4).
Let $\tau_i$ be the minimal face of  $\sigma_{J_i}$
containing $\sigma_{J_1} \cap \sigma_{J_2}$ for $i = 1, 2$.
There exists $\bmb \in \sigma_{J_1} \cap \sigma_{J_2} \cap \bQ^n$
such that $\bmb \in {\rm rel.int}(\tau_i)$ for $i = 1, 2$.
Put $J = J_{\bmb}(A)$.
Since $\tau_1, \tau_2 \in \Delta(A)$ by (a),
we know $\tau_1 = \tau_2 = \sigma_J$.
Since $J$ contains $J_1$ and $J_2$,
we have $\sigma_{J_1} \cap \sigma_{J_2} = \sigma_J$.

We have completed the proof of Theorem~\ref{FanStructure}.
\qed

\begin{Example}\label{NagataConj}
\begin{rm}
Let $\pi : X \rightarrow \bP_\bC^2$ be the blow-up at very general points
$p_1$, \ldots, $p_t$. 
Let $H$ be the hyperplane in $\bP_\bC^2$.
We put $E_i = \pi^{-1}(p_i)$ for $i = 1, 2, \ldots, t$.
Consider the multi-section ring (see Section~7)
\[
R = R(X; -E_1-E_2-\cdots-E_t, \pi^{-1}(H)) .
\]
Nagata~\cite{Nagata} conjectured that,
if $t \ge 10$ and $d \le \sqrt{t}m$, then $R_{(m,d)} = 0$.
If Nagata's conjecture is true, then
\[
\bRo(1, \sqrt{t}) + \bRo(0,1)
\]
is a maximal ray ideal cone of $R$.
Remark that, if $\sqrt{t} \not\in \bQ$, 
then it is not a rational polyhedral cone.
\end{rm}
\end{Example}

\section{Multi-section rings defined by $\bQ$-divisors}

Let $X$ be a normal projective variety over a field $k$ of $\dim X > 0$.
Let $K$ be the function field of $X$.
Let $H_1(X)$ be the set of closed subvarieties of $X$ of codimension $1$.
Let $D_i$ be a $\bQ$-Weil divisor for $i = 1, 2, \ldots, n$, that is a finite sum
\[
D_i = \sum_{F \in H_1(X)} m_{i, F}F ,
\]
where $m_{i, F} \in \bQ$.

We define a $\bZ^n$-graded ring as 
\begin{align*}
R(X; D_1, D_2, \ldots, D_n) &
= \bigoplus_{r_1, \ldots, r_n \in \bZ} H^0(X, {\mathcal O}_X(\sum_i r_iD_i))
t_1^{r_1}t_2^{r_2} \cdots t_n^{r_n} \\
& \subset K[t_1^{\pm 1}, t_2^{\pm 1}, \ldots, t_n^{\pm 1}] ,
\end{align*}
where
\[
H^0(X, {\mathcal O}_X(\sum_i r_iD_i))
= \{ a \in K^\times \mid {\rm div}(a) + \sum_i r_iD_i \ge 0 \} \cup \{ 0 \}
\subset K .
\]
We sometimes denote $R(X; D_1, D_2, \ldots, D_n)$ simply by $R$.
We put 
\[
m_{i, F} = \frac{p_{i, F}}{q_{i, F}} ,
\]
where $p_{i, F} \in \bZ$ and $q_{i, F} \in \bN$ such that ${\rm GCD}(p_{i, F}, q_{i, F}) = 1$.
Here, remark that $q_{i,F} = 1$ if $p_{i,F} = 0$.
We put
\[
q_F = {\rm LCM}( q_{1, F}, q_{2, F}, \ldots, q_{n, F} ) .
\]
We define
\[
P_F
= \bigoplus_{r_1, \ldots, r_n \in \bZ} H^0(X, {\mathcal O}_X(\sum_i r_iD_i - \frac{1}{q_F}F))
t_1^{r_1}t_2^{r_2} \cdots t_n^{r_n} 
\subset R(X; D_1, D_2, \ldots, D_n) .
\]

In this section, we shall prove the following:

\begin{Theorem}\label{EKW}
Let $X$ be a normal projective variety over a field $k$ with $\dim X > 0$.
Let $D_1$, $D_2$, \ldots, $D_n$ be $\bQ$-Weil divisors.
\begin{enumerate}
\item
The ring $R(X; D_1, D_2, \ldots, D_n)$ is a Krull domain.
The set of height $1$ homogeneous prime ideals of $R(X; D_1, D_2, \ldots, D_n)$
is contained in 
$\{ P_F \mid F \in H_1(X) \}$.
\item
Assume that $a_1 D_1 + a_2 D_2 + \cdots + a_n D_n$ is an ample Cartier
divisor for some $a_1, a_2, \ldots, a_n \in \bZ$.
Then the set of height $1$ homogeneous prime ideals of $R(X; D_1, D_2, \ldots, D_n)$
coincides with
$\{ P_F \mid F \in H_1(X) \}$.
Putting $\bma = (a_1, a_2, \ldots, a_n)$, the ray ideal $J_{\bma}(R(X; D_1, D_2, \ldots, D_n))$ is not contained in any height $1$ prime ideal of $R(X; D_1, D_2, \ldots, D_n)$.
\item
Assume that $a_1 D_1 + a_2 D_2 + \cdots + a_n D_n$ is an ample Cartier
divisor for some $a_1, a_2, \ldots, a_n \in \bZ$.
We put
\[
\{ F_1, F_2, \ldots, F_{\ell} \}
= \{ F \in H_1(X) \mid \mbox{$m_{i,F} \neq 0$ for some $i$} \} 
\]
and
\[
M  = \bigoplus_{j = 1}^\ell \frac{1}{q_{F_j}}\bZ \supset
L  = \sum_{i = 1}^n \bZ (m_{i, F_1}, m_{i, F_2}, \ldots, m_{i, F_\ell}) .
\]
Then we have an exact sequence
\begin{equation}\label{exact0}
L \cap \bZ^\ell \longrightarrow {\rm Cl}(X) \longrightarrow 
{\rm Cl}(R(X; D_1, D_2, \ldots, D_n)) \longrightarrow  \frac{M}{L + \bZ^\ell} \longrightarrow 0 ,
\end{equation}
where ${\rm Cl}( \ )$ denotes the divisor class group.
\end{enumerate}
\end{Theorem}

By the above theorem, we can immediately prove the following corollary.

\begin{Corollary}\label{UFD}
Assume that $a_1 D_1 + a_2 D_2 + \cdots + a_n D_n$ is an ample Cartier
divisor for some $a_1, a_2, \ldots, a_n \in \bZ$.
Then $R(X; D_1, D_2, \ldots, D_n)$ is factorial if and only if
\begin{itemize}
\item
$M = L + \bZ^\ell$, and
\item
${\displaystyle \left\{ \left. \sum_{i=1}^{n} b_i D_i \ \right|  
\mbox{
\begin{tabular}{l}
$b_{1}$, $b_{2}$, \ldots, $b_{n}$ are integers such that \\
$\sum_{i = 1}^{n}b_{i}(m_{i, F_1}, m_{i, F_2}, \ldots, m_{i, F_\ell}) \in \bZ^\ell$
\end{tabular}
}
\right\} }$ generates
$ {\rm Cl}(X)$.
\end{itemize}
\end{Corollary}

Now, we start to prove Theorem~\ref{EKW}.

First, we shall prove (1).
Let $v_{F}$ be the normalized valuation of the discrete valuation ring
${\mathcal O}_{X,F}$ for $F \in H_1(X)$.
We put
\begin{align*}
R_F & = \bigoplus_{r_1, \ldots, r_n \in \bZ}
\{ a \in K \mid v_F(a) + \sum_ir_im_{i,F} \ge 0 \} t_1^{r_1}t_2^{r_2}\cdots t_n^{r_n} \\
S & = K[t_1^{\pm 1}, t_2^{\pm 1}, \ldots, t_n^{\pm 1}] .
\end{align*}
By definition, we have
\begin{eqnarray}\label{riyuu1}
R =  \left( \bigcap_{F \in H_{1}(X)} R_{F} \right) \bigcap S .
\end{eqnarray}
Remark that $R_F$ satisfies $R \subset R_F
\subset S$.
Let $\alpha_F$ be a generator of the maximal ideal of the
discrete valuation ring ${\mathcal O}_{X,F}$.
Then we have
\[
\{ a \in K \mid v_F(a) + \sum_ir_im_{i,F} \ge 0 \}
= \alpha_F^{- \lfloor \sum_ir_im_{i,F} \rfloor } {\mathcal O}_{X,F} ,
\]
where $\lfloor \sum_ir_im_{i,F} \rfloor$ is the largest integer which is
not bigger than $\sum_ir_im_{i,F}$.

Here, we show that $R_F$ is a Noetherian normal domain.
If $m_{i,F}$ is an integer for all $i$ and $F$, then 
$R_F$ is a Noetherian normal domain since
\begin{align*}
R_F & = \bigoplus_{r_1, \ldots, r_n \in \bZ}
 \alpha_F^{- \sum_ir_im_{i,F}} {\mathcal O}_{X,F} t_1^{r_1}t_2^{r_2}\cdots t_n^{r_n} \\
 & = {\mathcal O}_{X,F}[ ( \alpha_F^{-m_{1,F}}t_1 )^{\pm 1} , ( \alpha_F^{-m_{2,F}} t_2)^{\pm 1}, \ldots,  ( \alpha_F^{-m_{n,F}}t_n )^{\pm 1}] .
\end{align*}
Suppose that $m_{i,F}$'s are rational numbers.
%We put
%\[
%q = {\rm LCM}(q_F \mid F \in H_1(X) ) .%%%
%\]
%Here remark that $q_{F} = 1$ except for finitely many $F$.
We put
%\begin{equation}\label{R_F^{(q_F)}}
\[
R_F^{(q_F)} = \bigoplus_{r_1, \ldots, r_n \in \bZ}
\{ a \in K \mid v_F(a) + \sum_iq_Fr_im_{i,F} \ge 0 \} t_1^{q_Fr_1}t_2^{q_Fr_2}\cdots t_n^{q_Fr_n} .
%\end{equation}
\]
Since $q_Fm_{i,F}$'s are integers, $R_F^{(q_F)}$ is a Noetherian normal domain.
For $a \in K$,
\begin{align*}
& \ \mbox{$a t_1^{r_1}t_2^{r_2}\cdots t_n^{r_n}$ is integral over $R_F^{(q_F)}$} \\
\Longleftrightarrow & \ \mbox{$(a t_1^{r_1}t_2^{r_2}\cdots t_n^{r_n})^{q_F}$ is integral over $R_F^{(q_F)}$} \\
\Longleftrightarrow & \ \mbox{$(a t_1^{r_1}t_2^{r_2}\cdots t_n^{r_n})^{q_F} \in R_F^{(q_F)}$} \\
\Longleftrightarrow & \ \mbox{$a t_1^{r_1}t_2^{r_2}\cdots t_n^{r_n} \in R_F$} .
\end{align*}
Therefore $R_{F}$ is the integral closure of $R_F^{(q_F)}$ in $S$.
Hence $R_F$ is a Noetherian normal domain.

Put
\[
Q_{F} = \bigoplus_{r_1, \ldots, r_n \in \bZ}
\{ a \in K \mid v_F(a) + \sum_ir_im_{i,F} > 0 \} t_1^{r_1}t_2^{r_2}\cdots t_n^{r_n} .
\]
Then $Q_{F}$ is a prime ideal of $R_{F}$ satisfying 
\begin{equation}\label{PF=QFcapR}
P_{F} = Q_{F} \cap R .
\end{equation}
It is easy to see $Q_{F} = \sqrt{\alpha_{F}R_{F}}$.
Remark that $Q_{F}$ is the unique height $1$ homogeneous prime ideal of $R_{F}$.
Furthermore, we have 
\begin{equation}\label{valu}
v_{Q_{F}}(\alpha_{F}) = q_{F} \ \ \mbox{and} \ \ v_{Q_F}(t_i) = q_Fm_{i,F}
\end{equation}
for $i = 1, 2, \ldots, n$, 
where $v_{Q_{F}}$ is the normalized valuation of the discrete valuation ring $(R_F)_{Q_{F}}$.
Let $\{ Q_{\lambda} \mid \lambda \in \Lambda \}$ be the set of non-homogeneous height $1$ prime ideals of $R_{F}$.
Since $R_{F}$ is a Krull domain,
\begin{equation}\label{riyuu2}
R_{F} = (R_{F})_{Q_{F}} \bigcap \left( \bigcap_{\lambda \in \Lambda} (R_{F})_{Q_{\lambda}} \right) .
\end{equation}
It is easy to see
\begin{equation}\label{riyuu3}
S = R_{F}[\alpha_{F}^{-1}] = \bigcap_{\lambda \in \Lambda} (R_{F})_{Q_{\lambda}} .
\end{equation}
Therefore we have
\begin{align*}
R = & \left( \bigcap_{F \in H_{1}(X)} R_{F} \right) \bigcap S = \left( \bigcap_{F \in H_{1}(X)} (R_{F})_{Q_{F}} \right) \bigcap S  \\ = & \left( \bigcap_{F \in H_{1}(X)} ((R_{F})_{Q_{F}} \cap Q(R)) \right) \bigcap (S \cap Q(R))
\end{align*}
by (\ref{riyuu1}), (\ref{riyuu2}) and (\ref{riyuu3}).
Since $R$ is the intersection of discrete valuation rings with finiteness condition (e.g., Section~12 in Matsumura~\cite{Mat}),
we know that $R$ is a Krull domain.
By Theorem~12.3 in \cite{Mat}, we know that the set of height one prime ideals of $R$ is contained in
\[
\{ P_F \mid F \in H_1(X) \} \cup \{ Q \cap R \mid Q \in {\rm Spec}(S), \ \height Q = 1 \} .
\]
The second assertion of (1) follows from the fact that any height $1$ prime ideal $Q$ of $S$ does not contain a non-zero homogeneous element of $R$. 

Next, we shall prove (2). 
Remark that $Q(R)$ coincides with $Q(S)$ since $\sum_{i}a_iD_i$ is ample.
We shall prove that $R_{P_F}$ coincides with $(R_F)_{Q_F}$ for any $F \in H_{1}(X)$.
We have
$R_{P_F} \subset (R_F)_{Q_F}$ by (\ref{PF=QFcapR}).
In order to prove the opposite inclusion, it is enough to show $R_F \subset R_{P_F}$.
Let $M_{F}$ be the set of homogeneous elements contained in $R \setminus P_{F}$.
We want to show $R_{F} \subset R[M_{F}^{-1}]$.
In the case where all the $m_{i,G}$'s are integers,
it is proved in the latter half of the proof of Theorem~1.1 in \cite{EKW}.
Here, put 
\[
q = {\rm LCM}\{ q_G \mid G \in H_1(X) \} .
\]
Then $R[M_{F}^{-1}]$ contains $R_{F}^{(q)}$,
where 
\[
R_F^{(q)} = \bigoplus_{r_1, \ldots, r_n \in \bZ}
\{ a \in K \mid v_F(a) + \sum_i q r_im_{i,F} \ge 0 \} t_1^{qr_1}t_2^{qr_2}\cdots t_n^{qr_n} .
\]
Since $R[M_{F}^{-1}]$ is integrally closed in $S$ and $R_{F}$ is integral over $R_{F}^{(q)}$,
$R[M_{F}^{-1}]$ contains $R_{F}$.
Thus we have 
\begin{equation}\label{7.1}
R_{P_F}=(R_F)_{Q_F} .
\end{equation}
Hence the set of height $1$ homogeneous prime ideal of $R(X; D_1, D_2, \ldots, D_n)$
coincides with
$\{ P_F \mid F \in H_1(X) \}$.
The latter assertion immediately follows from the first half.

Next, we shall prove (3). 
Let ${\rm Div}(X)$ be the set of Weil divisors on $X$, that is,
\[
{\rm Div}(X) = \left\{ \left. \sum_{F \in H_{1}(X)} n_{F} F \ \right| \ n_{F} \in  \bZ \right\} .
\]
Let ${\rm HDiv}(R)$ be the set of homogeneous Weil divisors of $R$, that is,
\[
{\rm HDiv}(R) = \left\{ \left. \sum_{F \in H_{1}(X)} n_{F} P_{F} \ \right| \ n_{F} \in  \bZ \right\} .
\]
Let 
\[
\phi : {\rm Div}(X) \longrightarrow {\rm HDiv}(R) 
\]
be a map defined by 
\begin{equation}\label{7.2}
\phi( F ) = q_{F} P_{F}  .
\end{equation}
Here, we have an exact sequence
\begin{equation}\label{exact1}
0 \longrightarrow {\rm Div}(X) \stackrel{\phi}{\longrightarrow} {\rm HDiv}(R)
\stackrel{\varphi}{\longrightarrow} M/\bZ^{\ell} \longrightarrow 0 ,
\end{equation}
where $\varphi$ is defined by $\varphi(\sum_{F}n_{F} F ) = 
\overline{(n_{F_{1}}/q_{F_{1}}, n_{F_{2}}/q_{F_{2}}, \ldots,
n_{F_{\ell}}/q_{F_{\ell}})}$.

For $a \in K^{\times}$ and $b \in Q(R)^{\times}$, we define
\[
{\rm div}_{X}(a) = \sum_{F \in H_{1}(X)} v_{F}(a) F , \ \ 
{\rm div}_{R}(b) = \sum_{P \in H_{1}(R)} v_{R_P}(b) P ,
\]
where $H_{1}(R)$ is the set of prime ideals of height $1$
and $v_{R_P}$ is the normalized valuation of the discrete valuation ring $R_P$.
Then we have
\begin{align*}
{\rm Cl}(X) & = {\rm Div}(X) / \{ {\rm div}_{X}(a) \mid a \in K^{\times} \} \\
{\rm Cl}(R) & = {\rm HDiv}(R) / \{ {\rm div}_{R}(b) \mid \mbox{$b$ is a non-zero homogeneous element of $S$} \} 
\end{align*}
by \cite{Sa}, \cite{KK}.
Here, remark that, for a non-zero homogeneous element $b$ of $S$,
$v_{R_P}(b) = 0$ if $P$ is a non-homogeneous prime ideal of $R$ of height one.
Put
\[
P(X) = \{ {\rm div}_{X}(a) \mid a \in K^{\times} \} .
\]
By (\ref{valu}), (\ref{7.1}) and (\ref{7.2}), we have 
\[
\phi({\rm div}_{X}(a)) = {\rm div}_{R}(a)
\]
 for any $a \in K^{\times}$.
Thus we have the following commutative diagram with exact rows: 
\begin{equation}\label{exact2}
\begin{array}{ccccccccc}
0 & \longrightarrow & {\rm Ker}(\xi) & \longrightarrow & \bZ^{n} & \stackrel{\xi}{\longrightarrow} & L & \longrightarrow & 0 \\
& & \phantom{\scriptstyle \eta'} \downarrow {\scriptstyle \eta'}   & & \phantom{\scriptstyle \eta} \downarrow {\scriptstyle \eta} & & \phantom{\scriptstyle \zeta} \downarrow {\scriptstyle \zeta} & \\
0 & \longrightarrow & {\rm Div}(X) / P(X) &
\stackrel{\overline{\phi}}{\longrightarrow} & {\rm HDiv}(R) / \phi(P(X)) & 
\stackrel{\overline{\varphi}}{\longrightarrow} & M/\bZ^{\ell} & \longrightarrow & 0 ,
\end{array}
\end{equation}
where the exact sequence in the second row is induced by (\ref{exact1}), 
$\zeta$ is induced by the inclusion $L \hookrightarrow M$,
$\eta$ is defined by 
\[
\eta(\bme_{i}) = {\rm div}_{R}(t_{i}) = \sum_{F \in H_{1}(X)} q_{F}m_{i,F} P_{F} ,
\]
$\xi$ is defined by 
\[
\xi(\bme_{i}) = (m_{i, F_1}, m_{i, F_2}, \ldots, m_{i, F_\ell}) .
\]
Then the induced map $\eta'$ coincides with $0$.
Since the cokernel of $\eta$ is ${\rm Cl}(R)$,
 we obtain the exact sequence as in (\ref{exact0}).
We have completed the proof of Theorem~\ref{EKW}.
\qed

\vspace{2mm}

Let us prove Corollary~\ref{UFD}.
The map
\[
L \cap \bZ^{\ell} \longrightarrow {\rm Cl}(X)
\]
in  (\ref{exact0}) sends 
\[
\sum_{i = 1}^{n}b_{i}(m_{i, F_1}, m_{i, F_2}, \ldots, m_{i, F_\ell}) \in L \cap \bZ^\ell
\ \ \ (b_{1}, b_{2}, \ldots, b_{n} \in \bZ)
\]
to 
\[
\sum_{i=1}^{n} b_i\left( \sum_{F \in H_{1}(X)}m_{i,F} F \right) = \sum_{i=1}^{n} b_i D_i .
\]
Therefore Corollary~\ref{UFD} immediately follows from Theorem~\ref{EKW} (3).
\qed

\section{Demazure construction for $\bZ^{n}$-graded Krull domains}

The aim of this section is to develop the Demazure construction
for  $\bZ^{n}$-graded Krull domains.

For a $\bZ^{n}$-graded domain $A$, $A_{(0)}$ denotes the homogeneous 
localization at $0$, that is, $A_{(0)} = A[M^{-1}]$, where
$M$ is the set consisting of all the non-zero homogeneous elements of $A$.
Remark that $A_{(0)}$ also has a structure of a $\bZ^{n}$-graded ring.

\begin{Lemma}\label{big}
Let $A$ be a $\bZ^{n}$-graded domain such that $A_{\bmo}$ is a field.
Take $\bma \in {\rm rel.int}(C(A)) \cap \bQ^{n}$.
Let $N$ be the set of non-zero homogeneous elements with degree
on the half line $\bRp   \bma$, that is,
\[
N = \left( \bigcup_{\bmb \in \bRp  \bma \cap \bZ^n} A_{\bmb} \right) \setminus \{ 0 \}  .
\]
Then $A[N^{-1}]$ coincides with $A_{(0)}$.
\end{Lemma}

\proof
Let $M$ be the set consisting of all the non-zero homogeneous elements of $A$.
For any $g \in M$, there exists $u \in A$ such that $gu$ is a homogeneous element contained in $N$
since $\bma$ is in the relative interior of $C(A)$.
Then the assertion immediately follows from $1/g = u/(gu)$.
\qed

\begin{Definition}\label{var}
\begin{rm}
Let $A$ be a $\bZ^{n}$-graded domain such that $A_{\bmo}$ is a field.
Take $\bma \in C(A) \cap \bQ^{n}$.
Assume that the ring
\[
A_{\bRo   \bma} = \bigoplus_{\bmb \in \bRo   \bma \cap \bZ^n} A_\bmb
\]
is Noetherian.
We think that $A_{\bRo   \bma}$ is a $\bN_0$-graded ring.
Then ${\rm Proj}(A_{\bRo   \bma})$ is called the {\em projective variety on the ray} $\bRo   \bma$,
and denoted by $X_{\bma}$.
\end{rm}
\end{Definition}

If $\bma \in {\rm rel.int}(C(A)) \cap \bQ^{n}$,
the function field of $X_{\bma}$ coincides with $(A_{(0)})_{\bmo}$ by Lemma~\ref{big}.
For any $\bmb \in C(A) \cap \bQ^{n}$,
the function field of $X_{\bmb}$ is contained in $(A_{(0)})_{\bmo}$.
Therefore, for $\bma \in {\rm rel.int}(C(A)) \cap \bQ^{n}$ and $\bmb \in C(A) \cap \bQ^{n}$,
there is a unique rational map from $X_{\bma}$ to $X_{\bmb}$ if the coordinate rings are Noetherian.

\begin{Lemma}\label{morph}
Let $A$ be a $\bZ^{n}$-graded domain such that $A_{\bmo}$ is a field.
Let $\bma$ and $\bmb$ be in $C(A) \cap \bQ^{n}$ such that
$A_{\bRo   \bma}$ and $A_{\bRo   \bmb}$ are Noetherian.
Assume that $J_{\bmb}(A)$ contains $J_{\bma}(A)$.
\begin{enumerate}
\item
The function field of $X_{\bmb}$ is naturally contained in that of $X_{\bma}$.
\item
The unique rational map from $X_{\bma}$ to $X_{\bmb}$ is a morphism of schemes.
\end{enumerate}
\end{Lemma}

\proof
Let $M$ be the set consisting of all the non-zero homogeneous elements of $A$.
Let $N_{\bma}$ (resp.\ $N_{\bmb}$) be the set consisting of all the non-zero homogeneous elements of $A$ with degree in $\bRp   \bma$ (resp.\ $\bRp   \bmb$).
Then the function field of $X_{\bma}$ (resp.\ $X_{\bmb}$) is $(A[N_{\bma}^{-1}])_{\bmo}$ 
(resp.\ $(A[N_{\bmb}^{-1}])_{\bmo}$).
Here, remark that both $(A[N_{\bma}^{-1}])_{\bmo}$ and $(A[N_{\bmb}^{-1}])_{\bmo}$ are contained in the field
$(A[M^{-1}])_{\bmo}$.

We may assume that $I_{\bma}(A)^{m} \subset I_{\bmb}(A)$ for some $m \in \bN$ since 
$A_{\bRo   \bma}$ is Noetherian.
Therefore there exist non-zero homogeneous elements $b_{1}$, $b_{2}$, \ldots, $b_{s}$ with degree in $\bRp   \bmb$,
and non-zero homogeneous elements $c_{1}$, $c_{2}$, \ldots, $c_{s}$ such that
\begin{itemize}
\item
$\deg(c_{i}b_{i})$ is in $\bRp   \bma$ for $i = 1, 2, \ldots, s$, and
\item
$c_{1}b_{1}$, $c_{2}b_{2}$, \ldots, $c_{s}b_{s}$ generate an ideal of $A$ containing $I_{\bma}(A)^{m}$.
\end{itemize}
Then $c_{1}b_{1}$, $c_{2}b_{2}$, \ldots, $c_{s}b_{s}$ generate an ideal 
of $A_{\bRo   \bma}$ which contains $(A_{\bRp   \bma})^m$.
Therefore we have an affine covering
\[
X_{\bma} = \bigcup_{i=1}^{s}{\rm Spec}\left( (A_{\bRo   \bma}[(c_{i}b_{i})^{-1}])_{\bmo}  \right) .
\]
Since
\[
(A_{\bRo   \bmb}[b_{i}^{-1}])_{\bmo} = (A[b_{i}^{-1}])_{\bmo} \subset (A[(c_{i}b_{i})^{-1}])_{\bmo} = (A_{\bRo   \bma}[(c_{i}b_{i})^{-1}])_{\bmo} ,
\]
we have
\[
(A[N_{\bmb}^{-1}])_{\bmo} = Q((A_{\bRo   \bmb}[b_{i}^{-1}])_{\bmo}) \subset
Q((A_{\bRo   \bma}[(c_{i}b_{i})^{-1}])_{\bmo}) = (A[N_{\bma}^{-1}])_{\bmo} .
\]
The assertion (1) has been proved.

We have a morphism
\[
{\rm Spec}\left( (A_{\bRo   \bma}[(c_{i}b_{i})^{-1}])_{\bmo}  \right) \longrightarrow
{\rm Spec}\left( (A_{\bRo   \bmb}[b_{i}^{-1}])_{\bmo}  \right) \subset X_{\bmb} .
\]
Patching together these morphisms, we obtain a morphism $X_{\bma} \rightarrow X_{\bmb}$.
\qed

\begin{Definition}\label{VarietyRayIdealCone}
\begin{rm}
Let $A$ be a $\bZ^{n}$-graded domain such that $A_{\bmo}$ is a field.
Let $\sigma$ be a ray ideal cone of $A$ such that $A_{\sigma}$ is Noetherian.
Let $\bma$ and $\bmb$ be in ${\rm rel.int}(\sigma) \cap \bQ^{n}$.
Then, by Lemma~\ref{morph}, we have $X_{\bma} = X_{\bmb}$.
We denote this variety by $X_{\sigma}$, and call it the {\em variety of
the ray ideal cone} $\sigma$.
\end{rm}
\end{Definition}

\begin{Example}\label{flip}
\begin{rm}
Let $k$ be a field.
\begin{enumerate}
\item
Let $A = k[x,y,z]$ be a $\bZ^2$-graded polynomial ring with 
$\deg(x) = (-1,0)$, $\deg(y) = (0,1)$, $\deg(z) = (1,1)$ and $A_\bmo = k$.
Then $J_{(0,1)}(A) \supsetneq J_{(-1,1)}(A)$.
However, the morphism $X_{(-1,1)} \rightarrow X_{(0,1)}$ is an isomorphism.
\item
Let $A = k[x,y,z,w,u]$ be a $\bZ^2$-graded polynomial ring with 
$\deg(x) = \deg(y) = (1,0)$, $\deg(z) = \deg(w) = (0,1)$, $\deg(u) = (1,1)$ and $A_\bmo = k$.
There exist just two maximal chambers
$\bRo   (1,0) + \bRo   (1,1)$ and $\bRo   (1,1) + \bRo   (0,1)$.
Put $\bma = (2,1)$, $\bmb = (1,1)$ and $\bmc = (1,2)$.
Since $J_{\bmb}(A) \supset J_{\bma}(A)$ and $J_{\bmb}(A) \supset J_{\bmc}(A)$,
we have birational morphisms
$X_\bma \rightarrow X_\bmb$ and $X_\bmc \rightarrow X_\bmb$.
We know $X_\bmb = {\rm Proj}(k[xz,xw,yz,yw,u])$.
Then $X_\bma \rightarrow X_\bmb$ is the blow-up along $V_+(xz,yz)$, and
$X_\bmc \rightarrow X_\bmb$ is the blow-up along $V_+(xz,xw)$.
\end{enumerate}
\end{rm}
\end{Example}

\vspace{2mm}

We shall prove the following theorem:

\begin{Theorem}\label{Demazure}
Let $A$ be a $\bZ^{n}$-graded domain such that $A_{\bmo}$ is a field.
Suppose that $\dim C(A) = n$.
Then the following three conditions are equivalent:
\begin{itemize}
\item[(\RMN{1})]
The ring $A$ is a Krull domain such that $(A_{(0)})_{\bme_{i}} \neq 0$ for $i = 1, 2, \ldots, n$.
Furthermore, there exists a chamber $\sigma$ of $A$ satisfying the following two conditions:
\begin{itemize}
\item[(a)]
The ray ideal $J_{\sigma}$ has height bigger than $1$.
\item[(b)]
The ring $A_{\sigma}$ is Noetherian.
\end{itemize}
\item[(\RMN{2})]
There exist a normal projective variety $X$ over $A_\bmo$,
$\bQ$-divisors $D_{1}$, $D_{2}$, \ldots, $D_{n}$ and
$\bmc = (c_{1}, c_{2}, \ldots, c_{n}) \in \bZ^n$
satisfying the following three conditions:
\begin{itemize}
\item[(a)]
$\sum_{i}c_{i}D_{i}$ is an ample Cartier divisor.
\item[(b)]
There exists a $n$-dimensional rational polyhedral cone $\rho$
such that $\bmc \in {\rm int}(\rho)$ and 
$A_\rho$ is Noetherian.
\item[(c)]
The ring $A$ is isomorphic to $R(X; D_{1}, D_{2}, \ldots, D_{n})$ as a 
$\bZ^{n}$-graded ring.
\end{itemize}
\item[(\RMN{3})]
There exist a normal projective variety $X$ over $A_\bmo$,
$\bQ$-divisors $D_{1}$, $D_{2}$, \ldots, $D_{n}$ and
$\bmc = (c_{1}, c_{2}, \ldots, c_{n}) \in \bZ^n$
satisfying the following three conditions:
\begin{itemize}
\item[(a)]
$\sum_{i}c_{i}D_{i}$ is an ample Cartier divisor.
\item[(b)]
There exists a positive integer $b$ such that $bD_{1}$, $bD_{2}$, \ldots, $bD_{n}$
are Cartier divisors.
\item[(c)]
The ring $A$ is isomorphic to $R(X; D_{1}, D_{2}, \ldots, D_{n})$ as a 
$\bZ^{n}$-graded ring.
\end{itemize}
\end{itemize}
\end{Theorem}

\begin{Remark}\label{DRem}
\begin{rm}
The variety $X$ is not determined uniquely in Theorem~\ref{Demazure}.
Let $A$ be the ring in Example~\ref{flip} (2).
Then $A$ is equal to the multi-section ring of $X_\bma$, $X_\bmb$ and $X_\bmc$.

In the case where $A$ is Noetherian,
obviously we may remove the conditions (b) in (I), and (b) in (II).

When we show ${\rm (I)} \Rightarrow {\rm (III)}$, we choose $\bmc$ which is contained in ${\rm int}(\sigma)$.

When we show ${\rm (II)} \Rightarrow {\rm (I)}$, we choose $\sigma$ such that $\sigma$ contains $\bmc$, but ${\rm int}(\sigma)$ may not contain $\bmc$.

When we show ${\rm (II)} \Rightarrow {\rm (III)}$, we sometimes need to change the variety $X$.

If $\dim X = 0$, we think that any $\bQ$-divisor is $0$, and $0$ is an ample Cartier divisor of $X$.

Using Lemma~\ref{big}, it is easy to see that the following four conditions are equivalent:
\begin{itemize}
\item
$A = A_{\bmo}[t_{1}^{\pm 1}, t_{2}^{\pm 1}, \ldots, t_{n}^{\pm 1}]$ with 
$\deg(t_{i}) = \bme_{i}$ for $i = 1, 2, \ldots, n$,
\item
the above condition (I) is satisfied with $J_{\sigma} = A$,
\item
the above condition (II) is satisfied with $\dim X = 0$,
\item
the above condition (III) is satisfied with $\dim X = 0$.
\end{itemize}
\end{rm}
\end{Remark}

Let us start to prove Theorem~\ref{Demazure}.
Assume the condition~(III).
We shall prove (II).
By (b) in (III), there exists 
an $n$-dimensional rational polyhedral simplicial cone $\rho$
such that $\bmc \in {\rm int}(\rho)$
and $\sum_ibc'_iD_i$ is ample Cartier divisor
for any $(c'_1,c'_2, \ldots, c'_n) \in \rho \cap \bZ^n$.
Then, by  Zariski's theorem~(Lemma~2.8 in Hu-Keel~\cite{HK}),
we know that $A_\rho$ is Noeherian.

Next, we shall prove ${\rm (II)} \Rightarrow {\rm (I)}$.
When $\dim X = 0$, this implication follows as in Remark~\ref{DRem}.
From now on, we assume $\dim X > 0$.
By Theorem~\ref{EKW} (1), $A$ is a Krull domain.
By (a) in (II), we know $(A_{(0)})_{\bme_{i}} \neq 0$ for $i = 1, 2, \ldots, n$ immediately.
Choose a chamber $\sigma$ of $A_\rho$ such that
\[
\bmc \in \sigma \subset \rho .
\]
Here, ${\rm int}(\sigma)$ need not contain $\bmc$.
Since $A_\sigma = (A_\rho)_\sigma$,  $A_\sigma$ is Noetherian by Remark~\ref{Noetherian} (3).
By Proposition~\ref{reduce} (1), 
$\sigma$ is also a chamber of $A$.
Suppose that a height one homogeneous prime ideal $P$ of $A$ contains $J_\sigma$.
By Theorem~\ref{EKW} (2), there exists $F \in H_1(X)$ such that
$P = P_F$.
Take $\bma = (a_1,a_2,\ldots,a_n) \in {\rm int}(\sigma) \cap \bZ^n$
such that $\sum_{i = 1}^na_im_{i,G}$ is an integer for any $G \in H_1(X)$.
Since $\bmc \in \sigma$ and $\bma \in {\rm int}(\sigma)$,
$\bma + m\bmc \in {\rm int}(\sigma)$ for any $m \in \bN$.
Therefore $A_{\bma + m\bmc} \subset P_F$.
On the other hand,
since $\sum_ic_iD_i$ is ample and $\sum_{i = 1}^na_iD_i$ is a Weil divisor with integer coefficients, it is easy to see that $A_{\bma + m\bmc} \not\subset P_F$ for $m \gg 0$.
Contradiction.

Next, we shall prove ${\rm (I)} \Rightarrow {\rm (III)}$.
We may assume $J_{\sigma} \neq A$ by Remark~\ref{DRem}.
Changing a coordinate, we may assume that the interior of the chamber $\sigma$ contains
the unit vectors $\bme_{1}$, $\bme_{2}$, \ldots, $\bme_{n}$.
By Lemma~\ref{morph}, we have
\[
X_{\bme_{1}} = X_{\bme_{2}} = \cdots = X_{\bme_{n}} .
\]
We denote this scheme by $X$.
Then $X$ is a normal projective variety of dimension positive
with function field $(A_{(0)})_\bmo$ by Lemma~\ref{big}.
(If $\dim X = 0$, then $A_{\bRo   \bme_{1}}$ is one dimensional $\bN_0$-graded normal domain.
It is easy to see that $A_{\bRo   \bme_{1}}$ is a polynomial ring over the field $A_{\bmo}$ with one variable.
Hence $I_{\bme_1}(A)$ is a principal ideal.
In this case, $J_{\bme_{1}}(A)$ is a unit ideal if the height of $J_{\bme_{1}}(A)$ is bigger than $1$.)
Here, we can choose $0 \neq t_{i} \in (A_{(0)})_{\bme_{i}}$ for $i = 1, 2, \ldots, n$ by (I).
Then, by the Demazure construction~\cite{D}, \cite{W}, there exists $\bQ$-divisors
$D_{1}$, $D_{2}$, \ldots, $D_{n}$ of $X$ such that
\begin{itemize}
\item
there exists a positive integer $b$ such that $bD_{1}$, $bD_{2}$, \ldots, $bD_{n}$
are ample Cartier divisors on $X$,
\item
for $x \in (A_{(0)})_{\bmo}$, $m \in \bN$ and $i = 1, 2, \ldots, n$, 
\[
xt_{i}^{m} \in A_{m\bme_{i}} \Longleftrightarrow x \in H^{0}(X, \oo_{X}(mD_{i})) .
\]
\end{itemize}

Both $A$ and $R(X; D_{1}, D_{2}, \ldots, D_{n})$ are $\bZ^{n}$-graded subrings of $A_{(0)} = (A_{(0)})_\bmo[t_1^{\pm 1}, t_2^{\pm 1}, \ldots, t_n^{\pm 1}]$.
It is enough to show $A = R(X; D_{1}, D_{2}, \ldots, D_{n})$ as a subset of $A_{(0)}$.
We denote $R(X; D_{1}, D_{2}, \ldots, D_{n})$ simply by $R$.

By Proposition~\ref{reduce} (3), $(\bRo)^{n}$ is a chamber of 
the Noetherian $\bZ^n$-graded ring $A_{(\bRo)^{n}}$.
Let $B$ be the subring of $A_{(\bRo)^{n}}$ generated by
\[
\bigcup_{i=1}^{n} \left[ \bigcup_{m > 0}A_{m \bme_{i}}\right]
\]
over $A_{\bmo}$.
Then the inclusion $B \rightarrow A_{(\bRo)^{n}}$ is finite by Theorem~\ref{Arai}.
Therefore $A_{(\bRo)^{n}}$ is the integral closure of $B$ in $A_{(0)}$.
On the other hand, $R$ is integrally closed in $A_{(0)}$ by Theorem~\ref{EKW},
and $R_{(\bRo)^{n}}$ is integral over $B$ 
by Zariski's theorem (Lemma~2.8 in Hu-Keel~\cite{HK}).
Therefore we have
\[
A_{(\bRo)^{n}} = R_{(\bRo)^{n}} .
\]
Let $y_{1}$, $y_{2}$, \ldots, $y_{t}$ be a homogeneous generating system of 
the homogeneous maximal ideal of $A_{\bRo  (\bme_{1}+ \bme_{2}+ \cdots+\bme_{n})}$.
Here, for a subring $T$ of $A_{(0)}$, we define the ideal transform to be
\[
D_{\underline{y}}(T) = 
\{ x \in A_{(0)} \mid \mbox{ $y_{1}^{m}x, \ y_{2}^{m}x, \ldots, \  y_{t}^{m}x \in T$ for $m \gg 0$ } \} .
\]
By definition, $D_{\underline{y}}(A_{(\bRo)^{n}})$ contains $A$.
Therefore we obtain $D_{\underline{y}}(A_{(\bRo)^{n}}) = D_{\underline{y}}(A)$.
In the same way, we obtain $D_{\underline{y}}(R_{(\bRo)^{n}}) = D_{\underline{y}}(R)$.
Thus we have
\[
D_{\underline{y}}(A) = D_{\underline{y}}(A_{(\bRo)^{n}}) =
D_{\underline{y}}(R_{(\bRo)^{n}}) = D_{\underline{y}}(R) .
\]
By (I) (a), the ideal $(y_{1}, y_{2}, \ldots, y_{t})A$ is of height bigger than $1$.
Since $A$ is a Krull domain, we obtain $D_{\underline{y}}(A) = A$.
By Theorem~\ref{EKW} (2), $(y_{1}, y_{2}, \ldots, y_{t})R$ is of height bigger than $1$.
Since $R$ is a Krull domain by Theorem~\ref{EKW} (1), we obtain $D_{\underline{y}}(R) = R$.

We have completed the proof of Theorem~\ref{Demazure}.

\begin{Example}\label{exDemazure}
\begin{rm}
Let $A = k[x,y,z,w]$ be a $\bZ^{2}$-graded polynomial ring 
over a field $A_\bmo = k$ with $\deg(x) = (1,0)$, $\deg(y) = (2,0)$,
$\deg(z) = (0,1)$ and
$\deg(w) = (0,2)$.
Let $\sigma$ be a $2$-dimensional cone contained in 
\[
\bRo   (1,0) + \bRo   (0,1) .
\]
Then $A$ and $\sigma$ satisfy the condition (I) in Theorem~\ref{Demazure}.
In this case, $X_{\sigma} = \bP_{k}^{1} \times \bP_{k}^{1}$.
Letting
\begin{align*}
D_{1} & = \frac{1}{2} \left( \{ (1:0) \} \times \bP_{k}^{1} \right) \ \ \mbox{and} \\
D_{1} & = \frac{1}{2} \left( \bP_{k}^{1}  \times \{ (1:0) \} \right) ,
\end{align*}
$A$ is isomorphic to $R(X_{\sigma};D_{1},D_{2})$.

Let $B = k[x,y,z,u,v]$ be a $\bZ^{2}$-graded polynomial ring over a field $k$ 
such that $\deg(x) = (1,0)$, $\deg(y) = (2,1)$,
$\deg(z) = (1,1)$,
$\deg(u) = (1,2)$ and
$\deg(v) = (0,1)$.
We put 
\[
\sigma_{1} = \bRo   (2,1) + \bRo   (1,1) , \ \ 
\sigma_{2} = \bRo   (1,1) + \bRo   (1,2) .
\]
Then both $\sigma_{1}$ and $\sigma_{2}$ satisfy the condition (I) in Theorem~\ref{Demazure}.
\end{rm}
\end{Example}

Next, we shall prove the following theorem:

\begin{Theorem}\label{Demazure2}
Let $A$ be a $\bZ^{n}$-graded domain such that $A_{\bmo}$ is a field.
Suppose that $\dim C(A) = n$.
Then the following three conditions are equivalent:
\begin{itemize}
\item[(I)]
The ring $A$ is a Krull domain such that $(A_{(0)})_{\bme_{i}} \neq 0$ for $i = 1, 2, \ldots, n$.
Furthermore there exists a chamber $\sigma$ of $A$ satisfying the following two conditions:
\begin{itemize}
\item[(a)]
For any height $1$ prime ideal $P$ containing $J_\sigma$,
$P$ contains $J_\bma(A)$ for any $\bma \in {\rm int}(C(A)) \cap \bQ^n$.
\item[(b)]
The ring $A_{\sigma}$ is Noetherian.
\end{itemize}
\item[(II)]
There exist a normal projective variety $X$ over $A_\bmo$,
$\bQ$-divisors $D_{1}$, $D_{2}$, \ldots, $D_{n}$ and
$\bmc = (c_{1}, c_{2}, \ldots, c_{n}) \in \bZ^n$
satisfying the following three conditions:
\begin{itemize}
\item[(a)]
$\bmc$ is in ${\rm int}(C(A))$ and
$\sum_{i}c_{i}D_{i}$ is an ample Cartier divisor.
\item[(b)]
There exists an $n$-dimensional rational polyhedral cone $\rho$
such that $\bmc \in {\rm int}(\rho)$ and 
$A_\rho$ is Noetherian.
\item[(c)]
There exists a rational polyhedral cone $\tau$ such that 
the ring $A$ is isomorphic to $R(X; D_{1}, D_{2}, \ldots, D_{n})_\tau$ as a 
$\bZ^{n}$-graded ring.
\end{itemize}
\item[(III)]
There exist a normal projective variety $X$ over $A_\bmo$,
$\bQ$-divisors $D_{1}$, $D_{2}$, \ldots, $D_{n}$ and
$\bmc = (c_{1}, c_{2}, \ldots, c_{n}) \in \bZ^n$
satisfying the following three conditions:
\begin{itemize}
\item[(a)]
$\bmc$ is in $C(A)$ and
$\sum_{i}c_{i}D_{i}$ is an ample Cartier divisor.
\item[(b)]
There exists a positive integer $b$ such that $bD_{1}$, $bD_{2}$, \ldots, $bD_{n}$
are Cartier divisors.
\item[(c)]
There exists a rational polyhedral cone $\tau$ such that 
the ring $A$ is isomorphic to $R(X; D_{1}, D_{2}, \ldots, D_{n})_\tau$ as a 
$\bZ^{n}$-graded ring.
\end{itemize}
\end{itemize}
\end{Theorem}

\proof
First, we shall prove ${\rm (III)} \Rightarrow {\rm (II)}$.
Assume that {\rm (III)} is satisfied.
Then it is easy to see that there exists $\bmc' = (c'_{1}, c'_{2}, \ldots, c'_{n}) \in \bZ^n$
such that $\bmc'$ is in ${\rm int}(C(A))$ and $\sum_{i = 1}^nc'_iD_i$ is an ample Cartier divisor.
Then this implication will be proved in the same way as in the proof of Theorem~\ref{Demazure}.

Next, we shall prove ${\rm (II)} \Rightarrow {\rm (I)}$.
Put $R = R(X; D_{1}, D_{2}, \ldots, D_{n})$.
By Theorem~\ref{Demazure}, $R$ satisfies the condition (I) in Theorem~\ref{Demazure}.
Therefore there exists a chamber $\sigma$ of $R$ such that the height of $J_{\sigma}(R)$ is bigger than $1$
and $R_{\sigma}$ is Noetherian.
We may assume $\bmc = (c_{1}, c_{2}, \ldots, c_{n}) \in \sigma$ (cf.\ Remark~\ref{DRem}).
By (II) (a), $\bmc$ is in ${\rm int}(C(A))$.
Replacing $\sigma$ by $\sigma \cap C(A)$, we may assume that 
\[
\sigma \subset C(A) \subset \tau .
\]
Since $A = R_{\tau}$, we have
\[
A = R \cap (R_{(0)})_{\tau} .
\]
In order to show that $A$ is a Krull domain, it is enough to show that $(R_{(0)})_{\tau}$ is a Krull domain.
(The intersection of two Krull domains are Krull again.)
Since $\tau$ is an $n$-dimensional rational polyhedral cone, there exist a finite number of vectors
$\bmb_{1}$, $\bmb_{2}$, \ldots, $\bmb_{m}$ in $\bQ^{n}$ such that
\[
\tau = \bigcap_{i=1}^{m}\{ \bma \in \bR^{n} \mid (\bma, \bmb_{i}) \ge 0 \} .
\]
We may assume that $\tau \cap \{ \bma \in \bR^{n} \mid (\bma, \bmb_{i}) = 0 \}$ is a face of $\tau$
of dimension $n-1$ for $i = 1, 2, \ldots, m$.
Here, we put
\[
T_{i} = \bigoplus_{\bma \in \bZ^{n}, \ (\bma, \bmb_{i}) \ge 0 } (R_{(0)})_\bma
\]
and 
\[
Q_{i} = \bigoplus_{\bma \in \bZ^{n}, \ (\bma, \bmb_{i}) > 0 } (R_{(0)})_\bma .
\]
Then $Q_{i}$ is a prime ideal of the ring $T_{i}$. 
Put $V_{i} = (T_{i})_{Q_{i}}$.
Then
\[
T_i = R_{(0)} \cap V_i \ \ \mbox{and} \ \ (R_{(0)})_{\tau} = R_{(0)} \cap \left( \bigcap_{i=1}^{m} T_{i} \right)
= R_{(0)} \cap \left( \bigcap_{i=1}^{m} V_{i} \right) .
\]
Since $R_{(0)}$ is Noetherian normal domain and $V_{i}$'s are discrete valuation rings,
$(R_{(0)})_{\tau}$ is a Krull domain.
Hence $A$ is a Krull domain.
Here, we have 
\[
A = R \cap \left( \bigcap_{i=1}^{m} V_{i} \right) .
\]
Let $P$ be a height one prime ideal of $A$ containing $J_{\sigma}$.
Height one prime ideals which come from $R$ does not contain $J_{\sigma}$ 
by the definition of $\sigma$.
Thus there exists $i$ such that $P = A \cap Q_{i}V_{i}$.
Then
\[
P = \bigoplus_{\bma \in \bZ^{n}, \ (\bma, \bmb_{i}) > 0 } A_\bma .
\]
Thus $P$ satisfies the condition (I) (a).

Next, we shall prove ${\rm (I)} \Rightarrow {\rm (III)}$.
If the height of $J_\sigma$ is bigger than $1$, this implication immediately follows from Theorem~\ref{Demazure}.

Assume hat $\height J_\sigma = 1$.
Let $P$ be a height one prime ideal containing $J_{\sigma}$.
Then, by  the condition (I) (a), 
$P$ contains $I_\bma(A)$ for any $\bma \in {\rm int}(C(A)) \cap \bQ^n$.
Since $P$ contains a non-zero homogeneous element and the height of $P$ is one,
$P$ is a homogeneous prime ideal.
Put $\nu = C(A/P)$.
By the assumption, we have 
$\nu \cap {\rm int}(C(A)) = \emptyset$. 
Then we obtain 
\begin{equation}\label{empty}
(\nu - \nu) \cap {\rm int}(C(A)) = \emptyset .
\end{equation}
We shall prove the following claim:

\begin{Claim}\label{rationalhypersurface}
\begin{enumerate}
\item
The dimension of  $\nu$ is $ n-1$.
\item
There exists $\bmo \neq \bmb \in \bQ^n$ such that
\[
A = \bigoplus_{\bma \in \bZ^n, \ (\bma,\bmb)\ge 0}A_\bma \ \ \ \mbox{and} \ \ \
P = \bigoplus_{\bma \in \bZ^n, \ (\bma,\bmb) > 0}A_\bma .
\]
\end{enumerate}
\end{Claim}

First, we shall prove (1).
Let $s$ be the dimension of $\nu$.
By (\ref{empty}), we have $s < n$.
Assume that $s < n-1$.
Changing the basis of $\bZ^n$, we may assume
\[
\nu - \nu = \bR \bme_1 +  \bR \bme_2 + \cdots + \bR \bme_s .
\]
Let
\[
\phi : \bR^n \longrightarrow \bR^{n-s}
\]
be the projection defined by 
\[
\phi((d_1, d_2, \ldots, d_n)) = (d_{s+1}, d_{s+2}, \ldots, d_n) .
\]
By (\ref{empty}), we know that $\phi(C(A)) \neq \bR^{n-s}$.
Then there exists a non-zero vector~\footnote{
The difficulty in this proof lies in that we have to find such a vector $\bmb'$ contained in $\bQ^{n-s}$.
}
$\bmb' \in \bR^{n-s}$
such that $(\bma', \bmb') \ge 0$ for any $\bma' \in \phi(C(A))$.
Here, put
\begin{align*}
& \Omega  = \{ \phi(\bma) \mid \bma \in \bZ^n, \ A_\bma \neq 0 \} \\
& \Omega_P  = \{ \phi(\bma) \mid \bma \in \bZ^n, \ 0 \neq A_\bma \subset P \} .
\end{align*}
If $A_{\bma} \not\subset P$, then $\bma \in \nu$ and $\phi(\bma) = \bmo$. 
Thus we have $\Omega = \Omega_P \cup \{ \bmo \}$.
Remark that 
\begin{equation}\label{OmegaP}
\mbox{$\Omega_P$ generates $\phi(C(A))$ as a cone}
\end{equation}
and
\[
\Omega_P \subset \Omega \subset \phi(C(A)) \subset 
\{ \bma' \in \bR^{n-s} \mid (\bma', \bmb') \ge 0 \} .
\]
If $\phi(\bma) \neq \bmo$, then $A_{\bma} \subset P$.
Thus we have 
\[
P \supset \bigoplus_{\bma \in \bZ^n, \ (\phi(\bma), \bmb')>0} A_\bma .
\]
Remark that the right-hand-side is a non-zero prime ideal of $A$.
Since the height of $P$ is one, we know that
$P$ coincides with the right-hand-side. 
That is, for $0 \neq q \in A_\bma$,
\begin{equation}\label{equiv}
q \in P \Longrightarrow (\phi(\bma), \bmb') > 0
\Longrightarrow \phi(\bma) \neq \bmo
\Longrightarrow q \in P .
\end{equation}

Let $x = x_1 + x_2 + \cdots + x_\ell$ be a non-zero element of $A$, where $0 \neq x_i \in A_{\bmc_i}$ for $i = 1, 2, \ldots, \ell$, and assume that $\bmc_1$, $\bmc_2$, \ldots, $\bmc_\ell$ are mutually distinct.
Here, we define
\begin{align*}
{\rm ord}(x) & = \min\{ (\phi(\bmc_i), \bmb') \mid i = 1, 2, \ldots, \ell \} \ge 0 \\
{\rm in}(x) & = \sum_{ (\phi(\bmc_i), \bmb') = {\rm ord}(x)} x_i .
\end{align*}
Then it is easy to see that, for non-zero elements $x, y \in A$, we have
\begin{align*}
{\rm ord}(x) + {\rm ord}(y) & = {\rm ord}(xy) , \\
{\rm in}(x) \cdot {\rm in}(y) & ={\rm in}(xy) .
\end{align*}
Let $\bmc'$ be an element in $\bZ^{n-s}$.
We say that $x$ is $\phi$-homogeneous with $\phideg(x) = \bmc'$ 
if $\phi(\bmc_i) = \bmc'$ for $i = 1, 2, \ldots, \ell$.

For $0 \neq q \in A$, we have
\begin{equation}\label{equiv2}
q \not\in P \Longleftrightarrow {\rm ord}(q) = 0
\Longleftrightarrow \mbox{${\rm in}(q)$ is $\phi$-homogeneous with
 $\phideg({\rm in}(q)) = \bmo$}
\end{equation}
by (\ref{equiv}).
 
Since $A_P$ is a discrete valuation ring,
there exists $\pi \in P$ such that $PA_P = \pi A_P$.
Take $0 \neq y \in A_\bma \subset P$.
Since $y \in PA_P = \pi A_P$,
there exists $r \in \bN$ and $z, w \in A \setminus P$ such that
\begin{equation}\label{zpi=xy}
z \pi^r = w y .
\end{equation}
Then we know that 
\begin{equation}\label{phi}
{\rm in}(z) \cdot {\rm in}(\pi)^r = {\rm in}(w) \cdot y .
\end{equation}
By (\ref{equiv2}), 
${\rm in}(z)$ and
${\rm in}(w)$ are $\phi$-homogeneous 
with $\phideg({\rm in}(z)) = \phideg({\rm in}(w)) = \bmo$.
Therefore ${\rm in}(\pi)^r$ is $\phi$-homogeneous with
$\phideg({\rm in}(\pi)^r) = \phi(\bma)$.
Then it is easy to see that
${\rm in}(\pi)$ is $\phi$-homogeneous and
\[
\phideg(y) = r \cdot \phideg({\rm in}(\pi)) .
\]
% the right-hand-side of (\ref{phi}) is  $\phi$-homogeneous with $\phideg({\rm in}(w) \cdot y) = \phi(\bma)$.
%Then, both ${\rm in}(z)$ and ${\rm in}(\pi)$ are  $\phi$-homogeneous such that
%\[
%\phideg({\rm in}(z)) + \phideg({\rm in}(\pi)) = \phi(\bma) .
%\]
%Let $v$ be a non-zero homogeneous component of ${\rm in}(z) \in A$.
%Suppose $v \in A_\bmd$.
%Then, we have
%\[
%\phi(\bmd) + \phideg({\rm in}(\pi)) = \phi(\bma) .
%\]
%If $v \not\in P$, we have $\phi(\bma) = \phideg({\rm in}(\pi))$ by (\ref{equiv2}).
%If $v \in P$, we repeat the same argument replacing $y$ by $v$.
%Here, remark that
%\[
%0 \le {\rm ord}(v) = {\rm ord}({\rm in}(z)) = {\rm ord}(z) = 
%{\rm ord}(y) - {\rm ord}(\pi) < {\rm ord}(y) .
%\]
%By induction on ${\rm ord}(y)$, we can prove $\phi(\bma)$ is
%equal to $r \cdot \phideg({\rm in}(\pi))$ for some $r \in \bN$.
Thus we have
\[
\Omega_P = \{ r \cdot \phideg({\rm in}(\pi)) \mid r \in \bN \} .
\]
It contradicts that $\dim \phi(C(A)) = n-s \ge 2$ (cf.\ (\ref{OmegaP})).

Next we shall prove (2) of Claim~\ref{rationalhypersurface}.
Since $s = n-1$, we may assume that the vector $\bmb'$ in the proof of (1)
is contained in $\bQ^1$.
Put $\bmb = (\bmo, \bmb') \in \bQ^n$.
Then we have $(\phi(\bma), \bmb') = (\bma, \bmb)$ for any $\bma \in \bR^n$.
Hence we have
\begin{align*}
& A = \bigoplus_{\bma \in \bZ^n, \ (\phi(\bma), \bmb')\ge 0}A_\bma =  \bigoplus_{\bma \in \bZ^n, \ (\bma,\bmb)\ge 0}A_\bma , \\ 
& P = \bigoplus_{\bma \in \bZ^n, \ (\phi(\bma), \bmb') >  0}A_\bma =  \bigoplus_{\bma \in \bZ^n, \ (\bma,\bmb) > 0}A_\bma .
\end{align*}

We have completed the proof of Claim~\ref{rationalhypersurface}.

\vspace{2mm}

Since $A$ is a Krull domain, there exists only finitely many height one prime ideals containing $J_\sigma$.
Let $\{ P_{1}, P_{2}, \ldots, P_{m} \}$ be the set of height one prime ideals of $A$ 
containing $J_\sigma$.
By Claim~\ref{rationalhypersurface} (2),
there exist 
$\bmb_{1}$, $\bmb_{2}$, \ldots, $\bmb_{m}$ in $\bQ^{n}$ such that
\begin{equation}\label{P_i}
A = \bigoplus_{\bma \in \bZ^n, \ (\bma,\bmb_i)\ge 0}A_\bma \ \ \mbox{and} \ \ 
P_{i} = \bigoplus_{\bma \in \bZ^n, \ (\bma,\bmb_i) > 0}A_\bma
\end{equation}
for $i = 1, 2, \ldots, m$.
We put 
\[
E_i = \bigoplus_{\bma \in \bZ^n, \ (\bma,\bmb_i)\ge 0}(A_{(0)})_\bma \ \ \mbox{and} \ \ 
U_i = \bigoplus_{\bma \in \bZ^n, \ (\bma,\bmb_i)> 0}(A_{(0)})_\bma .
\]
Since $A \subset E_i$ and $P_i = A \cap U_i$, we have $A_{P_i} \subset (E_i)_{U_i}$.
Since $A_{P_i}$ is a discrete valuation ring, we have $A_{P_i} = (E_i)_{U_i}$.
In particular, we have 
\[
A_{P_i} \supset E_i
\]
for each $i$.

Choose $\bma \in {\rm int}(\sigma) \cap \bQ^{n}$.
Take a homogeneous generating system $\{ y_{1}, y_{2}, \ldots, y_{t} \}$ of the ideal $I_{\bma}(A_\sigma)$ with degree in $\bRp   \bma$.
Consider the ideal transform $D_{\underline{y}}(A)$.
By definition, it is a $\bZ^{n}$-graded ring.
Let $H_1(A)$ be the set of height one prime ideals of $A$ .
Since $A$ is Krull, we have
\[
A = \bigcap_{P \in H_1(A)} A_{P} 
\]
and
\begin{equation}\label{ht=1}
D_{\underline{y}}(A) = \bigcap_{P \in H_1(A), \ P \not\supset (\underline{y})A} A_{P} .
\end{equation}
Then we obtain
\[
A = D_{\underline{y}}(A) \cap \left( \bigcap_{i=1}^{m} A_{P_{i}} \right)
= D_{\underline{y}}(A) \cap \left( \bigcap_{i=1}^{m} 
E_i \right) .
\]
Here we define $\tau$ to be
\begin{equation}\label{tau}
\tau = \bigcap_{i=1}^{m}\{ \bma \in \bR^{n} \mid (\bma, \bmb_{i}) \ge 0 \} .
\end{equation}
Then
\[
A = D_{\underline{y}}(A) \cap (A_{(0)})_{\tau} = D_{\underline{y}}(A)_{\tau} .
\]
Since 
\[
\sigma \subset C(A) \subset \tau ,
\]
we know $D_{\underline{y}}(A)_{\sigma} = A_{\sigma}$, and it is Noetherian by (I) (b).
By Proposition~\ref{reduce} (3), $\sigma$ is a chamber of $D_{\underline{y}}(A)$.
By (\ref{ht=1}), there is no height one prime ideal of $D_{\underline{y}}(A)$
containing $J_{\bma}(D_{\underline{y}}(A))$ (cf. Theorem12.3 in \cite{Mat}).
Therefore $D_{\underline{y}}(A)$ satisfies the condition (I) in Theorem~\ref{Demazure}.
Then, by Theorem~\ref{Demazure}, 
there exist a normal projective variety $X$ over $A_\bmo$ and 
$\bQ$-divisors $D_{1}$, $D_{2}$, \ldots, $D_{n}$ 
satisfying conditions (a), (b) in (III) in Theorem~\ref{Demazure}, and
\[
D_{\underline{y}}(A) = R(X;D_{1},D_{2}, \ldots, D_{n}) .
\]
Therefore we have
\[
A = R(X;D_{1},D_{2}, \ldots, D_{n})_{\tau} .
\]
By Remark~\ref{DRem}, the vector $(c_{1},c_{2}, \ldots,c_{n})$ is in ${\rm int}(\sigma)$.
Therefore $(c_{1},c_{2}, \ldots,c_{n})$ is in $C(A)$.

We have completed the proof of Theorem~\ref{Demazure2}.
\qed

\begin{Example}
\begin{rm}
Let $A = k[x,y,z]$ be a $\bZ^2$-graded polynomial ring
with $\deg(x) = (1,0)$,  $\deg(y) = (1,1)$ and $\deg(z) = (0,1)$.
Then $A$ does not satisfy the condition~(I) in Theorem~\ref{Demazure2}.
\end{rm}
\end{Example}

\begin{Remark}
\begin{rm}
Let $A$ be a $\bZ^n$-graded Noetherian normal domain such that
$A_\bmo$ is a field and $(A_{(0)})_{\bme_i} \neq 0$ for $i = 1, 2, \ldots, n$.
Even if all the ray ideals of chambers are of height less than $2$,
$A$ is isomorphic to a section ring $R(X;D_1,D_2,\ldots, D_n)$ as follows.

By the map
\[
\bZ^n \simeq \bZ^{n} \times \{ 0 \} \hookrightarrow \bZ^{n+1} ,
\]
we think that $A$ is a $\bZ^{n+1}$-graded ring.
Let $\bme_i$ be the $i$th unit vector in $\bZ^n$.
Consider the polynomial ring
\[
B = A[x_{1}, x_{2}, \ldots, x_{n+1}, y_{1}, y_{2}, \ldots, y_{n+1}]
\]
with $\deg(x_i) = \deg(y_i) = (\bme_i, 1)$ for $i = 1, 2, \ldots, n$
and $\deg(x_{n+1}) = \deg(y_{n+1}) = (-\bme_1 - \bme_2 - \cdots - \bme_n, 1)$.
Let $\bma$ be a point in $\bQ^{n+1}$ sufficiently near $(0, \ldots, 0,1)$.
Then $J_\bma(B)$ contains $x_1x_2\cdots x_{n+1}$ and $y_1y_2\cdots y_{n+1}$, and therefore 
the height of $J_\bma(B)$ is at least $2$.
Then there exists a chamber $\sigma$ of $B$ such that $B$ and $\sigma$ satisfy 
the condition (I) in Theorem~\ref{Demazure}.
Then there exist a normal projective variety $X$ of $\dim X > 0$ and 
$\bQ$-divisors $D_{1}$, $D_{2}$, \ldots, $D_{n+1}$ such that
$B$ is equal to $R(X; D_1, D_2, \ldots, D_{n+1})$.
Therefore $A$ coincides with $R(X; D_1, D_2, \ldots, D_{n})$.
In this case, $X$ does not inherit properties of the ring $A$ as in the 
following example. 

Let $A = k[x]$ be a graded polynomial ring over an algebraically closed field $k$ with $\deg(x) = 1$.
Let $X$ be a blow-up of a smooth projective variety over $k$ at a closed point.
Let $E$ be the exceptional divisor of $X$.
Then we have $A = R(X;E)$.
\end{rm}
\end{Remark}

\vspace{3mm}

\noindent
\begin{tabular}{l}
Yusuke Arai \\
IBM Japan Services Company Ltd. \\
Nihonbashi, Chuo-ku, \\
Tokyo Hakozakicho 19-21, Japan 
\end{tabular}

\vspace{2mm}

\noindent
\begin{tabular}{l}
Ayaka Echizenya \\
Department of Mathematics \\
Faculty of Science and Technology \\
Meiji University \\
Higashimita 1-1-1, Tama-ku \\
Kawasaki 214-8571, Japan 
\end{tabular}

\vspace{2mm}

\noindent
\begin{tabular}{l}
Kazuhiko Kurano \\
Department of Mathematics \\
Faculty of Science and Technology \\
Meiji University \\
Higashimita 1-1-1, Tama-ku \\
Kawasaki 214-8571, Japan \\
{\tt kurano@meiji.ac.jp} \\
{\tt http://www.isc.ac.jp/\~{}kurano}
\end{tabular}


\begin{thebibliography}{99}

\bibitem{AH}
{\sc K. Altmann and J. Hausen},
{\em Polyhedral divisors and algebraic torus actions},
Math.\ Ann.\ {\bf 334} (2006), 557--607.

\bibitem{BH}
{\sc F. Berchtold and J. Hausen}, 
{\em GIT-equivalence beyond the ample cone},
Michigan Math.\ J.\ {\bf 54} (2006), 483--515.

\bibitem{Cu}
{\sc S. D. Cutkosky}, 
\emph{Symbolic algebras of monomial primes}, J.
reine angew.\ Math.\ \textbf{416} (1991), 71--89.

\bibitem{D}
{\sc M.~Demazure},
{\em Anneaux gradu\'es normaux}, 
in Seminaire Demazure-Giraud-Teissier,
Singularites des surfaces, Ecole Polytechnique, 1979.
Introduction \`a la th\'eorie des singularit\'es, II,  35--68, 
Travaux en Cours, 37, Hermann, Paris, 1988. 

%\bibitem{CK}
%{\sc S D. Cutkosky and K. Kurano}, 
%\emph{Asymptotic regularity of powers of ideals of points in a weighted projective plane}, 
%Kyoto J. Math. {\bf 51} (2011), 25--45. 

\bibitem{Do}
{\sc I. V. Dolgachev}, 
{\em Automorphic forms and quasihomogeneous singularities}. 
Func. Anal.\ Appl. {\bf 9} (1975), 149–151.

\bibitem{EKW}
{\sc E. J. Elizondo, K. Kurano and K.-i. Watanabe},
{\em The total coordinate ring of a normal projective variety},
J. Algebra {\bf 276} (2004), 625--637.

%\bibitem{GK} {\sc J. L. Gonz\'alez and K. Karu},
%{\it Some non-finitely generated Cox rings},
%Compos. Math.  {\bf 152}  (2016),  984--996.

\bibitem{GNW}{\sc S. Goto, K. Nishida and K.-i. Watanabe},
{\it Non-Cohen-Macaulay symbolic blow-ups for space monomial curves and counterexamples to Cowsik's question},
Proc. Amer. Math. Soc. {\bf 120} (1994), 38--392.

%\bibitem{GW}{\sc S. Goto and K.-i. Watanabe},
%{\it On graded rings I},
%J. Math. Soc. Japan {\bf 30} (1978), 179--213.

\bibitem{HK} {\sc Y.~Hu and S.~Keel},
{\it Mori dream spaces and GIT},
Michigan Math J. {\bf 48} (2000), 331--348.

\bibitem{Hu} {\sc C. Huneke},
{\it Hilbert functions and symbolic powers},
Michigan Math. J. {\bf 34} (1987), 293--318.

\bibitem{HS} {\sc C. Huneke and I. Swanson},
{\it Integral Closure of Ideals, Rings, and Modules},
London Mathematical Society Lecture Note Series 336,
Cambridge University Press, 2006.

\bibitem{KK} {\sc Y.~Kamoi and K.~Kurano},
{\it On Chow groups of G-graded rings}, 
Comm. Alg. {\bf 31} (2003), 2141-2160.

\bibitem{K31} {\sc K.~Kurano}
{\it The divisor class groups and the graded canonical modules of multi-section rings}, 
Nagoya Math J. {\bf 212} (2013), 139--157.


\bibitem{LV} {\sc A.~Laface and M.~Velasco},
{\it A survey on Cox rings},
Geom Dedicata {\bf 139} (2009), 269--287.
 
\bibitem{Mat} {\sc H. Matsumura},
{\it Commutative ring theory},
Cambridge University Press, 1990.

%\bibitem{Mori} {\sc S.~Mori},
%{\it On a generalization of complete intersections},
%J. Math. Kyoto Univ. {\bf 15} (1975), 619--646.

\bibitem{Nagata}  {\sc M.~Nagata},
{\it On the 14-th Problem of Hilbert}, 
Amer. J. Math. {\bf 81} (1959), 766--772.

\bibitem{Okawa} {\sc S. Okawa},
{\it On images of Mori dream spaces},
Math. Ann. {\bf 364} (2016), 1315--1342.

\bibitem{P}
{\sc H. Pinkham},
{\em Normal surface singularities with $\bC^*$-action},
Math.\ Ann.\ {\bf 227} (1977), 183–193.

\bibitem{Sa} {\sc P. Samuel},
{\it Lectures on unique factorization domains}, 
Tata Inst. Fund. Res., Bombay, 1964.

\bibitem{W}
{\sc K.-i. Watanabe},
{\em Some remarks concerning Demazure's construction of normal graded rings},
Nagoya Math. J. {\bf 83} (1981),  203--211.

\end{thebibliography}
\end{document}